\documentclass[11pt,reqno]{amsart}
\usepackage{amssymb}

\usepackage[T1]{fontenc}
\usepackage[dvipsnames]{xcolor}

\usepackage{palatino}
\input amssym.def
\usepackage{amsmath, amsfonts}
\usepackage{amscd, mathtools}
\usepackage[mathscr]{eucal}
\usepackage{palatino}
\usepackage[colorlinks=true]{hyperref}
\hypersetup{citecolor=blue, linkcolor=blue}

\newfont{\cyrr}{wncyr10}

\newcommand{\N}{{\mathfrak N}}

\newcommand{\Q}{{\mathbb Q}}
\newcommand{\R}{{\mathbb R}}
\newcommand{\C}{{\mathbb C}}
\newcommand{\rO}{{\rm O}}

\newcommand{\K}{{\mathbf K}}
\renewcommand{\L}{{\mathbf L}}

\newcommand{\q}{{\mathfrak{q}}}
\newcommand{\h}{{\mathfrak h}}
\renewcommand{\d}{{\mathfrak d}}
\newcommand{\qq}{{\frac{1}{2}Q<\N\q\leq Q}}

\renewcommand{\mod}{{\, \rm mod \, }}
\newcommand{\thmref}[1]{Theorem~\ref{#1}}
\newtheorem{thm}{Theorem}[section]

\newtheorem{lem}[thm]{Lemma}
\newtheorem{cor}[thm]{Corollary}
\newtheorem{prop}[thm]{Proposition}

\newtheorem{rmk}{Remark}[section] 
\newtheorem{defn}{Definition}
\newtheorem{conj}{Conjecture}

\newcommand{\corref}[1]{Corollary~\ref{#1}}
\newcommand{\propref}[1]{Proposition~\ref{#1}}

\newcommand{\lemref}[1]{Lemma~\ref{#1}}

\renewcommand{\b}{{\mathfrak{b}}}
\renewcommand{\a}{{\mathfrak{a}}}
\renewcommand{\c}{{\mathfrak{c}}}
\newcommand{\e}{{\mathfrak{e}}}

\renewcommand{\P}{{\mathfrak{p}}}

\renewcommand{\O}{{{\mathcal{O}}}}

\newcommand{\p}{{\mathcal{P}}}

\renewcommand{\L}{\mathbf{L}}

\newcommand{\g}{{\mathfrak g}}

\begin{document}
	
\title[Euler-Kronecker constants]{ Explicit upper bounds on the average of \\ Euler-Kronecker constants
of narrow ray class fields}

\author{Neelam Kandhil, Rashi Lunia and
	Jyothsnaa Sivaraman}

\address{ Neelam Kandhil \\ \newline
	Max-Planck-Institut f\"ur Mathematik,
	Vivatsgasse 7, D-53111 Bonn, Germany.}

\address{ Rashi Lunia \\ \newline
	The Institute of Mathematical Sciences, A CI of Homi Bhabha National Institute, 
	CIT Campus, Taramani, Chennai 600 113, India.}

\address{Jyothsnaa Sivaraman \\ \newline
Statistics and Mathematics Unit,
Indian Statistical Institute, 
8th Mile, Mysore Rd, RVCE Post, Bengaluru, Karnataka 560059, India.}
\email{kandhil@mpim-bonn.mpg.de}
\email{rashisl@imsc.res.in}
\email{jyothsnaa.s@gmail.com}

\subjclass[2020]{11R18, 11R29, 11R42, 11Y60}
	% 11R18 Cyclotomic extensions, 11R29 Class numbers- class groups-discriminants, 11R42 Zeta functions and L-functions of number fields, 11Y60 Evaluation of number-theoretic constants

\keywords{Euler-Kronecker constants, Generalized Riemann Hypothesis,  Dedekind zeta functions, Narrow ray class fields, Generalized Dirichlet characters,
Hecke L-functions}

\begin{abstract} 
For a number field $\K$, the Euler-Kronecker constant $\gamma_\K$
associated to $\K$ is an arithmetic
invariant the size and nature of which is linked to some of the deepest
questions in number
theory. This theme was given impetus by Ihara who obtained bounds,
both unconditional as well as under GRH for Dedekind zeta functions. In this note, we study the analogous constants
associated to the narrow ray class fields of an imaginary quadratic field.
Our goal for studying such families is twofold.
First to show  that  for such families, the conditional bounds 
obtained by Ihara  can be improved on the average,
again under GRH for Dedekind zeta functions. Further, our  family of number fields are
non-abelian  while such average bounds have  earlier been studied
 for cyclotomic fields.  The  technical part of our work  is to
 make the dependence of these upper bounds on the ambient
number field explicit. Such explicit dependence is essential  to further our objective.

\end{abstract}

\maketitle   
\section{\large{Introduction}}
The Euler-Mascheroni constant $\gamma$, introduced by Euler \cite{LE}
 is given by  following limit
$$
\gamma  = \lim_{x \to \infty} \left( \sum_{n \le x } \frac{1}{n} - \log x \right).
$$
It is not known if $\gamma$ is rational
or irrational. It is perhaps prudent not to study  $\gamma$  in isolation,
but as a member of some  family whose behaviour can be studied en masse.
One such family was envisaged by Ihara \cite{YI2}  in 2005, the members
of which are referred to as Euler-Kronecker constants. To motivate 
these constants, let us note that $\gamma$ is  the
constant term in the Laurent series expansion of the Riemann
zeta function around $s=1$. Following this line of thought, 
one can define the Euler-Kronecker 
constant of a number field $\K$ as follows.

The Dedekind zeta function of a number
field $\K$ is given by 
$$
 \zeta_{\K}(s) = \sum_{\a \subseteq{\O_{\K}} \atop{ \a \neq (0)}} \frac{1}{\N\a^s}, \phantom{mm} \Re(s) > 1.
$$
It is known to have a meromorphic continuation to the entire complex plane
with only a simple pole at the point $s=1$. Thus we have a Laurent series
expansion for $\zeta_{\K}$ around $s=1$. This is given by
\begin{eqnarray}\label{Laurentseries}
\zeta_\K(s) =
\frac{\rho_\K}{s -1} + c_\K + \rO(s-1).
\end{eqnarray}
Following Ihara \cite{YI2}, we now define the Euler-Kronecker constant of $\K$ as
$$ \gamma_\K := c_\K/\rho_\K.$$
Since $\gamma_{\K}$ is evidently real, it is
natural to ask its sign. Ihara \cite{YI} 
conjectured that $\gamma_\K$ is positive when $\K$ is a cyclotomic field. 
However it has been shown by Ford, Luca and Moree \cite{FLM} that
$\gamma_{\Q(\zeta_p)}$ is negative for $p = 964477901$.
Further,  Ihara \cite{YI} showed that $\gamma_{\K}$ is negative infinitely
often under the Generalised Riemann Hypothesis (GRH)
for $\zeta_\K(s)$.

Another possible line of investigation is  to obtain bounds for $\gamma_{\K}$
in terms of other invariants of the ambient number field $\K$.
This was initiated  by Ihara in his papers  \cite{YI2}, \cite{YI}. 

For a number field $\K$ of degree  at least $3$, he showed 
 under GRH
 that
\begin{equation}\label{genbound}
-2(n_{\K}-1) \frac{ D_\K - n_{\K} + 1}{D_\K + n_{\K} -1} \left( \log \frac{ D_\K}{n_{\K}-1}+1 \right) -1 \leq \gamma_{\K} \leq   \left(\frac{D_\K +1}{D_\K  -1}\right)(2\log{D_\K}+1),
\end{equation}
where $D_\K=\log{\sqrt{|d_\K|}}$ and $n_{\K}$ and $d_{\K}$ denote the degree and discriminant
of $\K$ respectively. 

Recently, Dixit \cite{AD} showed the following unconditional bounds for $|\gamma_{\K}|$ for certain families of number fields.
We say that a number field $\K$ is almost normal
if there is a tower of number fields
$$
\K = \K_n \supset \K_{n-1} \supset \cdots \supset \K_1 = \Q, 
$$
such that $\K_{i+1}/\K_i$ is Galois for all $1 \le i \le n$.
\begin{thm}{\rm (Dixit, \cite{AD})}
Let $\K$ be an almost normal number field, not containing any quadratic subfields. Then
$$
\gamma_{\K} \ll (\log |d_{\K}|)^4 n_{\K}^3.
$$
On the other hand if $\K$ is a number field with solvable normal closure, not containing any
quadratic subfields. Then
$$
\gamma_{\K} \ll  (\log |d_{\K}|)^{c\log\log |d_{\K}|}
$$
for some absolute positive constant $c$.
\end{thm}

For an arbitrary number field, Dixit and Ram Murty \cite{DR}
related $\gamma_\K$ to the hypothetical
Siegel zero of $\zeta_{\K}(s)$. They showed the following.

\begin{thm}{\rm (Dixit and Ram Murty, \cite{DR,dixit-errata}}
For any number field $\K$, we have
$$
\gamma_{\K} \ll \log |d_{\K}|
$$
if $\zeta_{\K}(s)$ has no Siegel zero. If it has
a Siegel zero $\beta_0$, then the bound is
$$
\gamma_{\K} = \frac{1}{\beta_0(1-\beta_0)} + \rO(\log |d_{\K}|).
$$
\end{thm}

In the case of cyclotomic fields, Ihara suggested the following:
\begin{conj}{\rm(Ihara, \cite{YI2})}
There are positive constants $ 0 < a_0, a_1 \le 2$ such that for any $m$ sufficiently large and any $\epsilon >0$, we have
$$
(a_0 - \epsilon) \log m < \gamma_{\Q(\zeta_m)} < (a_1 + \epsilon ) \log m,
$$
where $\Q(\zeta_m)$ denotes the $m$-th
cyclotomic field.
\end{conj} 

\noindent The results of Ford, Luca and Moree \cite{FLM} indicate that $\gamma_{\Q(\zeta_m)}$ can be negative,
thus  the above lower bound  suggested by Ihara does not hold.

Assuming GRH for Dedekind zeta functions of cyclotomic fields, Ihara, Kumar Murty and Shimura \cite{IMS} proved that 
$$
\gamma_{\Q(\zeta_m)} \ll (\log m)^2.
$$

In 2010, Badzyan \cite{ABD} improved this result and showed that under GRH for Dedekind zeta functions of cyclotomic fields,
$$
\gamma_{\Q(\zeta_m)} \ll (\log m)(\log \log m).
$$

In 2011, Kumar Murty \cite{KM} proved that the upper bound of $|\gamma_{\Q(\zeta_p)}|$ coming from  Ihara's conjecture is true on average when restricted to primes. 
More precisely,

\begin{thm}{\rm (Kumar Murty, \cite{KM})}\label{avg}
	We have
	$$\frac{1}{\pi^{*}(Q)}\sum_{ \frac{1}{2}Q < p \le Q}|\gamma_{\Q(\zeta_p)}| \ll \log{Q},$$
	where the sum is over prime numbers $p$ in the interval $(\frac{1}{2}Q, Q] $ and $ \pi^*(Q)$ denotes the number of primes in this interval.
\end{thm}
Fouvry in \cite{Fou} studied the case of $\Q(\zeta_m)$ where $m$ is not necessarily prime and showed  the following.

\begin{thm}{\rm (Fouvry, \cite{Fou})}
Uniformly for $Q\geq 3$, we have
$$\frac{1}{Q}\sum_{\frac{1}{2}Q < m \le Q}\gamma_{\Q(\zeta_m)} = \log{Q}+\rO(\log{\log{Q}}).$$
\end{thm}

The aforementioned result of Fouvry has recently been refined by Hong, Ono and Zhang \cite{HKZ}
under the following conjecture.

\begin{conj}{\rm [Elliott and Halberstam conjecture \cite{EH}]~}\label{EHC}
For every real number $\theta < 1$ and for every positive 
integer $A> 0$, one has
$$
\sum_{q \leq x^\theta} 
 \max_{y \leq x} \max_{(a,q)=1} 
 \left| \sum_{n \le y \atop n \equiv a \bmod q} \Lambda(n) -   \frac{y}{\varphi(q)}\right| 
\ll_{A, \theta} \frac{x}{\log^A x}
$$
for all real numbers $x> 2$.
Here $\Lambda$ is the Von Mangoldt function.

\end{conj}

The precise result of Hong, Ono and Zhang is as follows.

\begin{thm}{\rm (Hong, Ono and Zhang, \cite{HKZ})}
Assume the Elliott-Halberstam conjecture, for $Q \to \infty$, we have
$$
\frac{1}{Q} \sum_{Q < m \le 2Q} |\gamma_{\Q(\zeta_m)} - \log m| = o(\log Q),
$$
where the sum is over integers $m$.
\end{thm}

We observe here that in a recently published work, Dixit and Ram Murty \cite{DR}
have given an alternative proof of the \thmref{avg} using certain new observations
for arbitrary number fields. Further they use the same method to show the following result.
\begin{thm}{\rm (Dixit and Ram Murty, \cite{DR})}
Assuming the Elliott-Halberstam conjecture, we have
$$
 \sum_{Q < p \le 2Q } |\gamma_{\Q(\zeta_p)} - \log p| = o(Q).
 $$
\end{thm}

Since the above results focus on cyclotomic fields, 
 the raison d'etre of our work
is to carry out such investigations   over families of 
Galois number fields which are  non-abelian.
One such family originates from
 narrow ray class fields of imaginary quadratic fields and it is the study of this family
 which we undertake. More precisely, we prove the following theorem. 
 
\begin{thm}\label{JR}
	Let $\K$ be an imaginary quadratic field with class number $h_{\K}$. For a non-zero principal prime ideal $\q$ of $\O_\K$, let $\K(\q)$ be the narrow ray class field modulo $\q$. 
	Assume   GRH for the Dedekind zeta functions of $\K(\q)$ for all non-zero principal prime ideals $\q$. Then for $ Q \geq 8 \exp(8 \cdot 10^{45} |d_{\K}|)$, we have
\begin{eqnarray}\label{finalbound}
\frac{1}{\pi^*(Q)}	\sideset{}{'}\sum_{\qq}|\gamma_{\K(\q)}|
 < |\gamma_{\K}| +
(6000h_{\K}^2+ 10^{17}h_\K  + 11) \log Q,
\end{eqnarray}
	where $'$ over the sum indicates that the sum is over principal prime ideals of $\O_\K$ and $ \pi^*(Q)$ denotes the number of
	principal prime ideals $\q$ of $\O_\K$ with norm in the interval $(\frac{1}{2}Q, Q] $.
\end{thm}

At this juncture, it is perhaps worthwhile to highlight a few points
about this result.
\begin{itemize}

\item The explicit dependence of the upper bound in \eqref{finalbound} on   $h_\K$ 
is what allows us to improve the lower bound in \eqref{genbound}  for infinitely many ray class fields of $\K$. 
To the best of our knowledge, this is the first instance of a non-abelian family over $\Q$ for which the bound has been improved. 

\item Even though our final result is conditional under GRH,  at various places, we have chosen
 unconditional bounds and not those under GRH unless essential.

\item We have not striven to optimise the numerical constants in the above theorem.
\end{itemize}

\bigskip
\noindent

This paper is organised as follows. In the next section we outline the strategy of our proof. 
In section 3, we list some results needed for our proofs 
while in  section 4 we derive some lemmas and propositions needed to prove our theorem.
The penultimate section gives a complete proof of \thmref{JR}.
Finally in the last section we show that there are infinitely many 
non-abelian fields in the family of narrow ray class fields considered by us.
\smallskip 
\section{Strategy of proof}
In this section, we quickly outline the strategy of the proof of \thmref{JR}.
In \cite{IMS}, Ihara, Kumar Murty and Shimura define, for a number field $\K$ and Hecke character $\chi$ modulo $\q$  of $\K$, (cf. section 3)
$$
\Phi_{\K, \chi} (x) = \frac{1}{x-1}\int_1^x \Big(\sum_{\a \subset \mathcal{O}_\K \atop 1 \leq \N\mathfrak{a} \leq t}  \frac{\Lambda(\mathfrak{a})}{\N\mathfrak{a}}\chi(\mathfrak{a})\Big)dt,~ \phantom{mm} (x> 1)
$$
and show that for $\chi \neq \chi_0$, ($\chi_0$ is the principal Hecke character)
$$
-\frac{L'}{L}(1, \chi) = \lim_{x \to \infty} \Phi_{\K, \chi} (x),
$$
where $L(s,\chi)$ is the Hecke L function corresponding to $\chi$.
We know that 
$$
\gamma_{\K(\q)} = \lim_{s \rightarrow 1}\left(
\frac{\zeta'_{\K(\q)}}{\zeta_{\K(\q)}}(s) + \frac{1}{s-1} \right) 
~~\text{ and }~~
\zeta_{\K(\q)}(s) = \zeta_{\K}(s) \prod_{\chi \neq \chi_0} L(s, \chi^*)
$$
where $\chi^*$ is the primitive character inducing the Hecke character $\chi$ modulo $\q$.
It now follows that
$$
|\gamma_{\K(\q)}| \le |\gamma_{\K}| +  \left|\sum_{\chi \neq \chi_0}\Phi_{\K,\chi}(x) \right| + \left|\sum_{\chi \neq \chi_0}\Big(\Phi_{\K,\chi}(x)-\Phi_{\K,\chi^*}(x)\Big) \right|+ \sum_{\chi \neq \chi_0} \left|\frac{L'}{L}(1,\chi^*)+\Phi_{\K, \chi^*}(x) \right|.
$$
The crux of the proof lies in estimating each of these sums on average as we vary norm of $\q$ between $Q/2$ and $Q$ where
$Q$ is a large real parameter. In \thmref{Phi1}, we estimate the average of the second term on the right by breaking the average into two parts. For small primes, we use an analogue of the Brun-Titchmarsh theorem (see \thmref{selbergsieve}) and some elementary estimates. 
For  large primes, we use the explicit Chebotarev density theorem under GRH due to Grenie and Molteni (see \thmref{explicitchebotarev}). The third term is relatively harmless and we estimate it directly (estimated after proof of \thmref{Phi1}).
To estimate the last term we use an explicit formula for 
$$
\frac{L'}{L}(1,\chi^*)+\Phi_{\chi^*}(x)
$$
proved by Ihara, Kumar Murty and Shimura in \cite{IMS} and deduce bounds under the assumption of 
GRH. \thmref{JR} is obtained by combining all these estimates.

\section{\large{Preliminaries}}\label{Prelims}

\smallskip 
\subsection{Notation}
Throughout this article, $\K$ will denote
an algebraic number field of degree $n_{\K}$, contained in $\C$. Let $\O_{\K}$ be its ring of integers and
$d_{\K}$ its discriminant with respect to $\Q$. Further, let $h_{\K}$ denote the class
number of $\O_{\K}$ and $\mu_\K$ the set of roots of unity in $\O_\K$. For any finite set $S$, we use $|S|$ to denote its cardinality. 

For an ideal $\a \subseteq \O_\K$, let $\N\a$  denote the order of the finite group $ \O_\K/\a$. 
We use $\p_{\K}$ to denote the set of
all prime elements of $\O_{\K}$. 
Further, $\P$
and $\q$ shall denote prime ideals in $\O_{\K}$.
We define the generalized Euler-phi function as follows: For any non-zero ideal $\a \subseteq \O_{\K}$,
$$
\varphi_{\K}(\a) = \N\a \prod_{\P \mid \a}\left( 1- \frac{1}{\N\P} \right).
$$
Let $\rho_{\K}$ denote the 
residue of the Dedekind zeta
function of $\K$ at $s=1$. For an imaginary quadratic field $\K$, the class number formula gives us
\begin{equation}\label{CNF}
\rho_\K=\frac{2\pi h_\K}{|\mu_\K|\sqrt{|d_\K|}}.
\end{equation}
The following result gives  bounds for $\rho_{\K}$ in terms of  $d_\K$.

\begin{lem}{\rm (Deshouillers, Gun, Ramar{\'e} and Sivaraman \cite{DGRS})}\label{DGRS}
We have
$$
 \frac{9}{25\sqrt{|d_{\K}|}}\leq
 \rho_{\K} \leq 6\left(\frac{2\pi^2}{5}\right)^{n_{\K}}|d_\K|^{1/4}.$$
\end{lem}

\smallskip
\noindent
We note that when $\K$ is imaginary quadratic,  the above can be refined
to obtain $$ \frac{\pi}{3\sqrt{|d_\K|}} \leq  \rho_{\K} \leq 6\left(\frac{2\pi^2}{5}\right)^2|d_\K|^{1/4}.$$

For number fields $\K \subset \L$, we recall that the relative discriminant $\Delta _{\L/\K}$ is defined as the $\O_\K$-ideal generated by the discriminants of all bases of $\L$ over $\K$ contained in $\O_{\L}$.  In this context, we record the following result.
\begin{lem}{\rm (Neukirch, page 202, \cite{JN})}\label{tower}
For a tower of fields $\Q \subset \K \subset \L$,	let $[ \L :\K ]$ denote the degree of $\L/\K$ and  $\Delta _{\L/\K}$ denote the relative discriminant of $\L$ over $\K$.  Then we have the following relation
$$
|d_{\L}| = |d_{\K}|^{[ \L ~:~ \K ]} ~ \N\Delta_{\L /\K}.
$$
\end{lem} 

Now we state the Conductor-Discriminant formula which will be used to bound the discriminant in \thmref{Phi1}.

\begin{thm} {\rm (Neukirch, page 534 (Conductor-Discriminant formula), \cite{JN})}\label{condisc}
 Let $\L/\K$ be a Galois extension of number fields and $\Delta _{\L/\K}$ denote its relative discriminant. Then we can express it by the following decomposition
 $$
\Delta_{\L/\K}
= \prod_{\chi ~~\textrm{irred.}}f(\chi)^{\chi(1)},
$$
where $\chi$ varies over the irreducible characters of the Galois group of $\L/\K$ and $f(\chi)$ denotes its Artin conductor (see page 527 of \cite{JN}). 
\end{thm}
Recall that when $\K = \Q$, $\zeta_\K(s)$ is the Riemann zeta function $\zeta(s).$ Let  $\frac{\Gamma'}{\Gamma}(s)$ and  $\frac{\zeta'}{\zeta}(s)$ denote the logarithmic derivatives of the gamma function and Riemann zeta function evaluated at $s$ respectively. Now we record the bounds on these functions.

\begin{lem}{\rm (Hall and Tenenbaum, page 146, \cite{hall})}\label{halllem}
For real $ \sigma > 1$, we have
$$
- \frac{\zeta'}{\zeta}(\sigma) < \frac{1}{\sigma - 1}.
$$
\end{lem}

\begin{lem}{ \rm (Ahn and Kwon, \cite{ahn})} \label{ahnlem}
Assume that $\Re(s) > \frac{1}{2}$. We have
$$ 
\Re \left(\frac{\Gamma'}{\Gamma}(s)\right) \leq 1.08 \log(|s| + 2).
$$
\end{lem}
For an arithmetic function $f$ and a positive arithmetic function $g$, $f(z) = \rO^*(g(z))$ means
that $|f(z)| \leq g(z)$.
 
\subsection{Narrow ray class groups and Hecke L-functions.}

Let $\c$ be a non-zero integral ideal of $\K$.
For any element $\alpha \in \K$, we say that  $\alpha ~\equiv~ 1\bmod^*\c$ if $\alpha$ satisfies the following
two conditions:
\begin{enumerate}
 \item The element $\alpha$ is of the form $a/b$
 for $a, b \in \O_{\K}$ where the ideals
 $(a)$ and $(b)$ are coprime to $\c$ and
 $
 a \equiv b   \bmod  \c.
 $

\item For any real embedding $\sigma$ of $\K$, $\sigma(\alpha)>0$.
\end{enumerate}
\begin{defn}
Let $I(\c)$ denote the set of non-zero fractional ideals relatively prime to $\c$ and $P_{\c}$ the group of non-zero principal fractional ideals $(\alpha) \subset \O_{\K}$ such that $\alpha \equiv  1 \bmod^* \c$.
 Then the group $I(\c)/P_{\c}$  is called  the narrow ray class group of $\K$ modulo $\c$, denoted  by $H_{\c}(\K)$.
\end{defn}

\begin{defn}
A character $\chi$, of $H_{\c}(\K)$ will be
called a generalized Dirichlet character
modulo $\c$. We shall use $\chi_0$ to
denote the principal generalized Dirichlet character.
\end{defn}
There is a natural restriction
from $H_{\c}(\K)$ to $H_{\c'}(\K)$
for all $\c' \mid \c$.
\begin{defn}
	We say that a character of 
	$H_{\c}(\K)$ is primitive if it
	does not factor through
	$H_{\c'}(\K)$ for every proper divisor $\c'$ of $\c$. Further, the conductor of
	a generalized Dirichlet character $\chi$
	modulo $\c$ is the smallest divisor
	$\c'$ of $\c$ such that $\chi$
	factors through a character
	modulo $\c'$.
\end{defn}

Further for a generalized
Dirichlet character $\chi$ modulo $\c$, we define the associated Hecke $L$- function as
$$
L(s, \chi, \K) = \sum_{0 \neq \a \subseteq{ \O_{\K}}\atop{(\a, \c)=\O_\K}} \frac{\chi([\a])}{\N\a^s}, \qquad \Re(s) > 1,
$$
where $[\a]$ denotes the class of the ideal $\a$ in $H_{\c}(\K)$. We will henceforth use $L(s, \chi)$ to denote $L(s, \chi, \K)$.
It is known that when $\chi$ is not the principal character, $L(s, \chi)$
has an analytic continuation to the entire complex plane.
If $\chi$ is principal then $L(s,\chi)$ has
a meromorphic continuation to the entire
complex plane with only a simple pole at $s=1$. The Generalized Riemann Hypothesis (GRH) for Dedekind zeta functions states that the
zeros of a Dedekind zeta function in the strip $0~<~\Re(s)~<~1$  lie on the line $\Re(s) = 1/2$.
Now we state
the following lemma of Lagarias and Odlyzko.
	\begin{lem}{\rm (Lagarias and Odlyzko, \cite{LO})}\label{diskolem}
	 Let $\chi$ be a generalized
		Dirichlet character modulo integral ideal $\c$
		and let $L(s, \chi)$ be the corresponding
		Hecke $L$-function. Then we have
		$$
		\frac{L'}{L}(s, \chi) + \frac{L'}{L}(s, \overline{\chi}) = \sum_\rho \left(\frac{1}{s - \rho} + \frac{1}{s- \overline{\rho}}\right) - \log ( |d_\K| \N\c') -  \mathbf{1}_{\chi}\left(\frac{2}{s} + \frac{2}{s-1}\right) - 2 ~ \frac{\Gamma'_\chi}{\Gamma_\chi}(s),
		$$
		where the sum is over non-trivial zeros of $L(s,\chi),$
		$ \mathbf{1}_{\chi}$ is $1$ if $\chi$ is a principal character otherwise it is $0$, $\c'$ is the conductor of $\chi$
		and for some fixed integer $ 0 \leq a_{\chi} \leq n_\K$,
		\begin{equation}\label{eq100}
			\Gamma_\chi(s) := 
			\Big[ \pi^{-\frac{s+1}{2}} \Gamma\left(\frac{s+1}{2}\right) \Big]^{a_{\chi}} ~ ~ 
			\Big[ \pi^{-\frac{s}{2}} \Gamma\left(\frac{s}{2}\right) \Big]^{n_\K - a_{\chi}}.
		\end{equation}
	\end{lem}
 Let $\chi \neq \chi_0$ be a  generalized Dirichlet character on an imaginary quadratic field $\K$ and $L(s, \chi)$ be the corresponding Hecke $L$- function. 
Now we consider the following function
\begin{equation} \label{eq101}
	\Phi_\chi(x):=\frac{1}{x-1}\int_1^x \Big(\sum_{\a \subset \mathcal{O}_\K \atop 1 \leq \N\mathfrak{a} \leq t}  \frac{\Lambda(\mathfrak{a})}{\N\mathfrak{a}}\chi([\mathfrak{a}])\Big)dt,
	\end{equation}
as defined in \cite{IMS}. Then we have the following 
result proved by Ihara, Kumar Murty and Shimura.
\begin{thm} {\rm (Ihara, Kumar Murty and Shimura, \cite{IMS})}\label{ims}
	For $x>1$ and any primitive generalized Dirichlet character $\chi$, we have
	\begin{equation*}
		\frac{L'}{L}(1, \chi) + \Phi_{\chi}(x)= \frac{1}{x-1}\sum_{\rho}\frac{x^{\rho}-1}{\rho(1-\rho)}+\log{\frac{x}{x-1}}+\frac{1}{x-1}\log{x},
	\end{equation*} 
	where the sum is over all non-trivial zeros of $L(s, \chi)$ counted with multiplicities, and
	\begin{equation*}
		\sum_{\rho}\frac{x^{\rho}-1}{\rho(1-\rho)}=\lim_{T \to \infty}\sum_{|\rm{Im}(\rho)|<T}\frac{x^{\rho}-1}{\rho(1-\rho)}.
	\end{equation*}
\end{thm}

\subsection{Counting integral ideals}

We now state a result on counting the number of integral ideals of $\O_{\K}$ for number field $\K$ due to 
Gun, Ramar{\'e} and Sivaraman. 
For any embedding $\sigma$ of $\K$, the Minkowski embedding
$\theta$ of $\K$ to $\R^2$ maps $x$ to  $(\Re(\sigma(x)), \Im(\sigma(x)))$.
In this setup, we have the following counting theorem.

\begin{thm}{\rm (Gun, Ramar{\'e} and Sivaraman, \cite{GRS2})}\label{counting}
Let $\a, \q$ be co-prime ideals of $\O_{\K}$,
$\mathfrak{C}$ be the ideal class of $\a\q$ in the class group of $\O_{\K}$
and $\Lambda(\a\q)$ be the lattice $\theta(\a\q)$ in $\R^2$, where $\theta$ 
is as defined above.  Also let
$$
S_{\beta}\left(\a, \q,  t^{2} \right)  
= \{\alpha \in \a~: |\theta(\alpha)|^2 \le t^{2}, \alpha \equiv \beta \bmod \q\}
$$
for some fix $\beta \in \O_{\K}$.
Then for any real number $t \ge 1$, we have
\begin{equation}\label{countingresult}
\left| S_{\beta}\left(\a, \q,  t^2 \right)  \right| 
 = 
\frac{2 \pi }{\sqrt{|d_{\K}|} \N(\a\q)} t^2 + ~
\rO^*\left( \frac{ 10^{13.66}\mathfrak{N}(\mathfrak{C}^{-1}) }
{ |\N(\a\q)|^{\frac{1}{2}}} t   +  1 \right),
\end{equation}
where 
$$
\N(\mathfrak{C}^{-1}) = 
\text{max}_{\b \in \mathfrak{C}^{-1}} 
\frac{1}{|\mathfrak{N}(\b)|^{\frac12}}.
$$
One can ignore $1$ in the error term when $\q = \O_{\K}$.
\end{thm}
As a corollary the authors of \cite{GRS2, GRS} deduce the following two theorems.
\begin{thm}{\rm (Gun, Ramar{\'e} and Sivaraman, \cite{GRS2})}\label{asymfinal}
Let $\K$ be an imaginary quadratic field.
Let $\q$ be an integral ideal of $\K$ and $[\b]$ be an element of $H_{\q}(\K)$.
For any real number $x \ge 1$, we have
\begin{equation*}
\sum_{\substack{\a\subset\O_\K\\ [\a]=[\b]\\ \N\a \le x}} 1 = \frac{\rho_{\K} \varphi(\q)}{|H_{\q}(\K)|}
\frac{x}{\N\q} 
+\rO^*\left(10^{21}\left(\frac{x}{\N\q}\right)^{1/2}
+  4\cdot 10^5 \right).
\end{equation*}
\end{thm}

\begin{thm}{\rm (Gun, Ramar{\'e} and Sivaraman, \cite{GRS})} \label{GRS}
Let $\K$ be an imaginary quadratic field. For any real number $x \geq 1$, 
$$\sum_{\a \subseteq \O_\K \atop \N\a \leq x} 1= \rho_\K x+\rO^*\left(10^{15}(h_\K\log(3h_\K))^{1/2}x^{1/2}\right).$$
\end{thm}

We now state a result of Garcia and Lee \cite{GLee} which will be used in the proof of \thmref{JR}.
\begin{thm}\label{GLee}
Let $\K$ be an imaginary quadratic field and $x \geq 2$. We have 
$$\sum_{\N\P\leq x}\frac{\log\N\P}{\N\P} =\log x +\rO^*\left(3+\frac{e^{75}|d_\K|^{1/3}(\log|d_\K|)^2}{\rho_\K} \right).$$
\end{thm}

Let $\L/\K$ be an abelian extension with Galois group $G$ and let $n_{\L}$ and $n_\K$ be the respective degrees over $\Q$. 

\begin{defn}
	Let $\P$ be a non-zero prime ideal of $\O_{\K}$,
	unramified in $\L$ and $\q$ be a prime ideal above it in $\O_{\L}$.
	The unique Galois element
	$\sigma$ in the Galois group of $\L/\K$ such that 
	$$
	\sigma(a) \equiv a^{\mathfrak{N}(\P)} \bmod \q ~\text{ for all } a \in \O_{\L}
	$$
	 is called the Artin symbol
	corresponding to $\P$ and the extension
	$\L/\K$. 
	The Artin symbol is extended multiplicatively to the set of all non-zero integral ideals of $\O_{\K}$ which are supported on the set of prime ideals which do not ramify in $\L/\K$.
	For such an integral ideal $\a$, we use $\left(\L/\K, \a\right)$ to
	denote the Artin symbol.
\end{defn}

Under the above notation, we have the following 
theorem from \cite{GM}.
\begin{thm}\label{explicitchebotarev}{\rm(Greni{\'e} and Molteni,  \cite{GM})}
Let $\sigma$ be an element of the group $G$ and let $\mathbf{1}_{\sigma}$ denote the characteristic function of $\sigma$. We define
$$
\psi(x, \L/\K, \sigma) = \sum_{\a \subset \O_{\K} \atop{  \P~~ \text{unramified}~~ \forall~~ \P\mid\a \atop{1 \le \N\a \le x}}} 
\mathbf{1}_{\sigma}((\a, \L/\K)) \Lambda(\a).
$$
Under GRH, for all $x \geq 1$,
\begin{equation*}
\Big||G|\psi(x, \L/\K, \sigma)-x \Big|\leq \sqrt{x} \left( \left(\frac{\log x }{2 \pi} +2\right) \log |d_{\L}| +
\left(\frac{\log^2 x }{8 \pi} +2\right) n_{\L} \right).
\end{equation*}
\end{thm}

The above result is an extension of their previous result for number fields \cite{GM2} which we shall now state
in a weaker form as required by our proof.

\begin{thm}\label{GM2}{\rm(Greni{\'e} and Molteni,  \cite{GM2})}
Let $\K$ be a number field. We define
$ \displaystyle
\psi(x) = \sum_{\a \subset \O_{\K},\atop{1 \le \N\a \le x}} 
 \Lambda(\a).
$
Under GRH, for all $x \geq 3$,
\begin{equation*}
|\psi(x)-x |\leq  7\log{|d_\K|}\sqrt{x}\log x+\sqrt{x}\log^2{x}+19\sqrt{x}.
\end{equation*}

\end{thm}

Before we proceed further, we would like to introduce the notion of the
narrow ray class field modulo a non-zero integral ideal $\c$ of $\O_{\K}$.
To do so, we first define the Artin map
in the following manner :
Given an abelian extension $\L/\K$ and a non-zero integral ideal $\c$ such that every prime of $\K$ that ramifies in $\L$ divides $\c$, we define $\Psi_{\L/\K, \c}  : 	I(\c)   \to  \operatorname{Gal}(\L/\K)$ as the map
\begin{eqnarray*}
 \displaystyle	\prod_{\P \nmid \c \atop{\P \text{ prime ideal of } \O_{\K}}} \P^{m_{\P}}  \mapsto   \prod_{\P \nmid \c \atop{\P \text{ prime ideal of } \O_{\K}}}(\L/\K, \P)^{m_{\P}}, \phantom{mm} m_{\P} \in \mathbb{N}
\end{eqnarray*}
where the integer $m_{\P}$ is $0$ for all but finitely
many primes $\P$.
It is known from class field theory
that this map is surjective.

\begin{defn}
For a non-zero ideal
$\c$, the unique extension $\L$
of $\K$ such that
\begin{enumerate}
	\item every prime ideal of $\K$ that ramifies in $\L$ divides $\c$ and 
	\item  the kernel of 
	$\Psi_{\L/\K, \c}$ is $P(\c)$
\end{enumerate}
is called the narrow ray class field of $\K$
with respect to $\c$. We shall
denote this field by $\K(\c)$.
\end{defn}

It is known that every finite abelian
extension $\L$ of an imaginary quadratic field $\K$ is contained in a 
narrow ray class field $\K(\c)$ of $\K$ for some integral ideal $\c$.
The greatest common divisor  of all such integral ideals
is called the conductor of $\L$ with respect to $\K$.

Corollary 1.2 of \cite{GM} restates the above
bounds in terms of the prime counting
function $\pi(x, \L/\K, \sigma)$, where
$$
\pi(x, \L/\K, \sigma) = \sum_{(0) \neq \P \subset \O_{\K}, \atop{ \P~~ \text{ unramified } \atop{1 \le \N\P \le x}}} 
\mathbf{1}_{\sigma}((\P, \L/\K)).
$$
In this setup, we now deduce a corollary
about principal prime ideals.
We first recall that the Hilbert class field of an imaginary quadratic field $\K$ 
is the narrow ray class field of $\K$ corresponding to the ideal $\O_{\K}$. By class field theory
it is also the unique abelian extension of
$\K$ where all the non-zero principal prime ideals
of $\O_{\K}$ split completely. 
For a number field $\K$, we use $\K(\O_{\K})$ to denote its Hilbert class field. 

\begin{cor}\label{principalprimecounting}
Let $\K$ be an imaginary quadratic field and $\K(\O_\K)$ the Hilbert class field of $\K$.
Let $\sigma_0$ be the trivial element of the Galois group of $\K(\O_\K)/\K$.
Under GRH, $\forall x \ge 2$,
\begin{eqnarray*}
\pi(x, \K(\O_\K)/\K, \sigma_0) & = & \Big|\{ \P \subseteq \O_{\K} : 1 \le  \N\P \le x, \P \text{ splits completely in}~ \K(\O_\K) \}\Big| \\
 & = &  \frac{1}{h_{\K}} \int_2^x \frac{1}{\log t} dt  + \rO^* \left( 5 \sqrt{x} \log |d_{\K}| + 2 \sqrt{x} \left(\frac{\log x}{8\pi} + 9 \right) \right). 
\end{eqnarray*}
\end{cor}

\smallskip

\section{\large{Requisite Propositions and Lemmas}}

In this section, we state and prove some lemmas and propositions required in the proofs of our theorems. The first lemma gives us a bound on the cardinality of $H_{\q}(\K)$ for  an imaginary quadratic field $\K$.
\begin{lem}\label{size}
	Let $\K$ be an imaginary quadratic field. For any non-zero prime ideal $\q$ in $\O_{\K}$,
	$$\frac{1}{12}h_\K  \N\q \leq |H_{\q}(\K)| < h_\K  \N\q.$$ 
\end{lem}
\begin{proof}
See Theorem 1 of  
 Chapter VI \cite{SL} (page 127).
\end{proof}

We now state and prove some lemmas to compute the sum of the inverse 
of the function $\varphi_{\K}$ introduced earlier.
For a number field $\K$ and a non-zero ideal $\a \subseteq \O_{\K}$,
we define  $$\displaystyle \sigma(\a) := \sum_{\b \subseteq \O_\K \atop \b \mid \a} \N\b.$$

\begin{lem}\label{lem10}
Let $\K$ be an imaginary quadratic field. Then for $x \geq 1$, we have
$$
\sum_{\a \subseteq \O_{\K} \atop{1 \le \N\a \le x}} \frac{\sigma(\a)}{\N\a} \leq \zeta_\K(2)\rho_\K x+ 10^{15}\zeta_{\K}(3/2) (h_\K\log(3h_\K))^{1/2} x^{1/2}.
$$
\end{lem}

\begin{proof}
We have
\begin{eqnarray*}
\sum_{\a \subseteq \O_{\K} \atop{1 \le \N\a \le x}} \frac{\sigma(\a)}{\N\a}
& = &\sum_{\a \subseteq \O_{\K} \atop{1 \le \N\a \le x}}  \frac{1}{\N\a}  \sum_{\b \subseteq \O_{\K} \atop{\b \mid \a}} \N\b  
=\sum_{\a \subseteq \O_{\K} \atop{1 \le \N\a \le x}}    \sum_{\b \subseteq \O_{\K} \atop{\b \mid \a}} \frac{1}{\N\b}
= \sum_{\b \subseteq \O_{\K} \atop{1 \le \N\b \le x}}  \frac{1}{\N\b} \sum_{\c \subseteq \O_{\K}  \atop{1 \le \N\c \le \frac{x}{\N\b}}} 1\\
 & \leq & \sum_{\b \subseteq \O_{\K} \atop{1 \le \N\b \le x}}  \frac{1}{\N\b}\left( \rho_\K \frac{x}{\N\b}+10^{15}(h_\K\log(3h_\K))^{1/2}\left(\frac{x}{\N\b}\right)^{1/2}\right)\\
 &\leq & \zeta_\K(2)\rho_{\K}x+ 10^{15}\zeta_{\K}(3/2)(h_\K\log(3h_\K))^{1/2} x^{1/2}.
\end{eqnarray*}
Here we have used \thmref{GRS} in the penultimate step. Using equation \eqref{CNF}, we get the desired result.

\end{proof}

\begin{lem}\label{lem11}
For a number field $\K$ and $ \a \subset \O_{\K}$, we have
$$\frac{\N\a}{\varphi_{\K}(\a)} \le \zeta_{\K}(2) \frac{\sigma(\a)}{\N\a}.$$
\end{lem}

\begin{proof}
For $\P \mid \a$, let $\alpha_{\P}$ be such that
$\P^{\alpha_{\P}} \| \a$.
We now consider
$$
\frac{\N\a}{\varphi_{\K}(\a)} = \prod_{\P \mid \a} \left(1 - \frac{1}{\N\P} \right)^{-1} = \prod_{\P \mid \a} \frac{1 + \frac{1}{\N\P}}{1 - \frac{1}{\N\P^2}} 	=\prod_{\P \mid \a}  \frac{\N\P^{\alpha_{\P}} + \N\P^{\alpha_{\P}-1} }{\N\P^{\alpha_{\P}} \left(1 - \frac{1}{\N\P^2}\right)} 
\le \zeta_{\K}(2) \frac{\sigma(\a)}{\N\a}.$$
\end{proof}

\begin{lem}\label{lem12}
Let $\K$ be an imaginary quadratic field. We now have for $x \geq 3,$
$$
\sum_{ \a \subset \O_{\K} \atop{1 \le \N\a \le x}} \frac{1}{\varphi_{\K}(\a)} 	
\leq \zeta_\K(2)^2 \rho_\K \log ex + 2 \cdot 10^{15}\zeta_{\K}(3/2)\zeta_\K(2)(h_\K\log(3h_\K))^{1/2}.
$$
\end{lem}

\begin{proof}
Consider
\begin{eqnarray*}
\sum_{\a \subseteq \O_{\K}  \atop{1 \le \N\a \le x}} \frac{1}{\varphi_{\K}(\a)} 
&=& \sum_{1 \le t \le x}  \sum_{ \a \subseteq \O_{\K}  \atop{ \N\a = t}} \frac{1}{\varphi_{\K}(\a)}
= \sum_{1 \le t \le x} \frac{1}{t}  \sum_{ \a \subseteq \O_{\K}  \atop{ \N\a = t}} \frac{\N\a}{\varphi_{\K}(\a)}.
\end{eqnarray*}
Let $\displaystyle a_t =  \sum_{ \a \subseteq \O_{\K}  \atop{ \N\a = t}} \frac{\N\a}{\varphi_{\K}(\a)}$ and $\displaystyle A(x)=\sum_{1 \le t \le x} a_t$. We now apply
partial summation formula to get
\begin{eqnarray}\label{partialsumn}
\sum_{\a \subseteq \O_{\K} \atop{1 \le \N\a \le x }} \frac{1}{\varphi_{\K}(\a)} 
&=& \sum_{1 \le t \le x} \frac{a_t}{t} = 	\frac{A(x)}{x} + \int_{1}^x \frac{A(u)}{u^2} du.
\end{eqnarray}
By  \lemref{lem10} and \lemref{lem11},
\begin{align*}
A(x)&= \sum_{1 \le t \le x} a_t = 
\sum_{1 \le t \le x} \sum_{ \a \subseteq \O_{\K}  \atop{ \N\a = t}} \frac{\N\a}{\varphi_{\K}(\a)} 
= \sum_{\a \subseteq \O_{\K} \atop{1 \le \N\a \le x }}  \frac{\N\a}{\varphi_{\K}(\a)}  \\
&\leq  \zeta_\K(2)^2 \rho_\K x+10^{15}\zeta_{\K}(3/2)\zeta_\K(2)(h_\K\log(3h_\K))^{1/2} x^{1/2}.
\end{align*}
Using the above bound in \eqref{partialsumn}, we get the desired result.

\end{proof}

Now we compute the upper bound of the absolute discriminant of the narrow ray class field $\K(\q)$ corresponding to modulus $\q$ of number field $\K$.

\begin{lem} \label{Disc}
Let $\K$ be an imaginary quadratic field. 	Let $\q$ be a non-zero prime ideal in $\O_{\K}$ and $d_{\K(\q)} $ be the discriminant of $\K(\q)$. Then,
$$\log |d_{\K(\q)}| \leq \rho_{\K}\sqrt{|d_\K|} \N\q \log(|d_{\K}|\N\q).$$
\end{lem}

\begin{proof}
Using \lemref{tower}, we have
$$
|d_{\K(\q)}|= |d_{\K}|^{|H_{\q}(\K)|}\N\Delta_{\K(\q)/\K}.
$$
 Since $\operatorname{Gal}(\K(\q)/\K)$ is abelian, every irreducible character of $\operatorname{Gal}(\K(\q)/\K)$ is of degree $1$.
If $\chi$ is a Hecke character of $\operatorname{Gal}(\K(\q)/\K)$, $f(\chi)$ is also the conductor of a subfield of $\K(\q)$ which contains $\K$ ~~(see page 535, Prop 11.10 of \cite{JN}). Therefore, $f(\chi)$ divides $\q$. By \thmref{condisc}, we have
$$\Delta_{\K(\q)/\K}
= \prod_{\chi ~~\textrm{irred.}}f(\chi) .$$
Taking norm on both
 sides, we get
$$\N\Delta_{\K(\q)/\K}=\prod_{\chi ~~\textrm{irred.}}\N f(\chi) \leq \prod_{\chi ~~\textrm{irred.}}\N\q= \N\q^{|H_{\q}(K)|}.$$
Using \lemref{size}, we get
\begin{align*}
\log |d_{\K(\q)}| =
|H_{\q}(\K)|\log{|d_{\K}|}+ \log \N\Delta_{\K(\q)/\K}
\leq \rho_\K\sqrt{|d_\K|} \N\q \log(|d_{\K}|\N\q).
\end{align*} 	

\end{proof}

\subsection{Selberg's sieve}

Let $(t)$ be a non-zero principal integral ideal of $\O_\K$. For a non-zero integral ideal $\b$, let $\b_{(t)}$ denote the
maximal divisor of $\b$ co-prime to $(t)$, that is
$$
\b_{(t)}=\prod_{(\P,(t))=1}\P^{\nu_\P(\b)},
$$
where $\nu_\P(\b)$ denotes the valuation of $\b$ at $\P$.
Also let $d(\b)$ denote the number of divisors of $\b$ and $\Omega_{\K}(\b)$ denote the number of prime ideals dividing $\b$ counted with multiplicity. For a non-zero integral ideal $\b$, we define 
$$
\rho(\b) = 2^{\Omega_{\K}(\b_{(t)})},
\phantom{mm}
f_{(t)}(\b) = \frac{\N\b}{2^{\Omega_{\K}(\b_{(t)})}}
\phantom{mm}  {\textrm{and}}
\phantom{m}
f_{1,(t)}(\b) = \sum_{\a \mid \b} \mu(\a) f_{(t)}\left( \frac{\b}{\a} \right).
$$
 Evidently all these functions are multiplicative.
Let 
$$
\mathbb{P}(z) = \prod_{ N\P \le z \atop{(\P, (2310))=1}} \P.
$$
Further, for an ideal $\e$ of $\O_{\K}$ coprime to $(2)$,  we define
\begin{eqnarray*}
S_{\e}(z) 
=
\sum_{ \a \subseteq \O_\K \atop{{1 \le \N(\a) \le z \atop{(\a, (2310)\e) = \O_{K}}}}} \frac{\mu^2(\a)}{f_{1,(t)}(\a)},
\phantom{m}\text{ and }\phantom{m}
\lambda_{\e} 
=
 \mu(\e) \frac{f_{(t)}(\e)S_{\e}(\frac{z}{\N(\e)})}{f_{1,(t)}(\e)S_{\O_{\K}}(z)}. 
\end{eqnarray*}
We observe here that $2310$ is the product of all the primes
upto $12$. Working with primes greater than $12$ allows us to get explicit lower bounds for certain sums
which appear in the Selberg sieve (see equation \eqref{Sok} below).
For non-zero integral ideals $\b_1$ and $\b_2$, we use $[\b_1, \b_2]$ to denote the 
least common multiple of the ideals $\b_1$
and $\b_2$. In this setup, we have the following result.

\begin{lem}\label{errorterm}
We have
$$
\sum_{\b_1, \b_2 \mid \mathbb{P}(z), \atop{ \N\b_i \le z}} 
|\lambda_{\b_1}\lambda_{\b_2}| \frac{{\rho([\b_1,\b_2])}}{\sqrt{\N([\b_1,\b_2])}}
~\le~
\zeta_{\K}\left(\frac32 \right)^{16} z.
$$
\end{lem}
\begin{proof}
The proof of this lemma follows along the same lines as Lemma 11
of \cite{GRS2}.
We consider the sum
$$
\sum_{\b_1, \b_2 \mid \mathbb{P}(z), \atop{ \N\b_i \le z}} |\lambda_{\b_1}
\lambda_{\b_2}| \frac{\rho([\b_1,\b_2])}{\sqrt{\N([\b_1,\b_2])}}
~=~
\sum_{\d \mid \mathbb{P}(z),  \atop{ \N\d \le z}}\frac{\sqrt{\N\d}}{\rho(\d)} 
\sum_{ \b_i \mid \mathbb{P}(z), \atop{\d = (\b_1,\b_2)
\atop{\N\b_i \le z}}} \frac{|\lambda_{\b_1} \lambda_{\b_2}| \rho(\b_1) 
\rho(\b_2)}{\sqrt{\N(\b_1\b_2)}}.
$$
From the expression of $\lambda_{\b}$ and with $y=z/\N\d$, we get
\begin{align*}
\sum_{ \substack{\N\c\le y,\\ (\c,\d (2310))=\O_\K}} \frac{|\lambda_{\d\c}| \rho(\c)}{\sqrt{\N\c}} 
&  \le 
  \frac{\sqrt{y}~\N\d}{\rho(\d)f_{1,(t)}(\d)}\prod_{\N\P > 11}\biggl(1+\frac{\rho(\P)}{\sqrt{\N\P}(\N\P-\rho(\P))}\biggr).
\end{align*}

We thus get 
\begin{align*}
  \sum_{\b_1, \b_2 \mid \mathbb{P}(z), \atop{ \N\b_i \le z}} |\lambda_{\b_1}
  \lambda_{\b_2}| \frac{\rho([\b_1,\b_2])}{\sqrt{\N([\b_1,\b_2])}}
&\le~
  z \prod_{\N\P > 11}\biggl(1+\frac{\rho(\P)}{\sqrt{\N\P}(\N\P-\rho(\P))}\biggr)^2
  \biggl(1+\frac{\rho(\P)\sqrt{\N\P}}{(\N\P-\rho(\P))^2}\biggr).
\end{align*}
Note that
$$
\prod_{\N\P > 11} \left(1+\frac{\rho(\P)}{\sqrt{\N\P}(\N\P-\rho(\P))} \right) 
~\le~
\zeta_{\K} \left( \frac32 \right)^{4}
\text{ and }
\prod_{\N\P > 11} \left(1+\frac{\rho(\P)\sqrt{\N\P}}{(\N\P-\rho(\P))^2} \right) 
~\le~
 \zeta_{\K}\left(\frac32 \right)^{8}.
$$

This completes the proof of the lemma.
\end{proof}

\begin{prop}\label{selbergsieve}
Let $\K$ be an imaginary quadratic field and $t$ be a fixed non-zero element in $\O_\K$, which is not a unit in $\O_{\K}$.
We have for 
$
u \ge  2 \exp\left(8 \cdot 10^{45}  \cdot |d_{\K}|\right)
$

\begin{eqnarray*}
\sum_{1 \le \N(\alpha) \le u \atop{\alpha, t\alpha+1~\in~ \p_{\K}}} 1 
 & \le &
\left(  \frac{2 \cdot 10^{10}}{ \sqrt{|d_{\K}|}   \rho_{\K}^2  } 
  \right) \frac{\N(t) u}{\varphi((t)) \log^2 \frac{ u^{1/4}}{2}}
  +
  10^{30} u^{3/4}.
\end{eqnarray*}
\end{prop}

\begin{proof}
Let $z$ be a real number such that $\sqrt{z} > 12$.
We would like to estimate the sum
\begin{eqnarray*}
\sum_{1 \le \N(\alpha) \le u \atop{\alpha, t\alpha+1 \in \p_{\K}}} 1 
& \le &
 \sum_{1 \le \N(\alpha) \le z \atop{\alpha, t\alpha+1 \in \p_{\K}}} 1 
  +
 \sum_{z < \N(\alpha) \le u \atop{1 \leq \N(t\alpha+1) \le z \atop{\alpha, t\alpha+1 \in \p_{\K}}}} 1 
 +
 \sum_{z < \N(\alpha) \le u \atop{((\alpha(t\alpha+1)), \mathbb{P}(z))=1}} 1 \\
& \le &
 \sum_{1 \le \N(\alpha) \le u \atop{((\alpha(t\alpha+1)), \mathbb{P}(z))=1}} 1  
 +
 4|\mu_{\K}|z.
\end{eqnarray*}
Let us consider the first sum
$$
 \sum_{ 1 \le \N(\alpha) \le u \atop{((\alpha(t\alpha+1)), \mathbb{P}(z))=1}} 1 
 =
  \sum_{1 \le \N(\alpha) \le u } \left( \sum_{\b \mid ((\alpha(t\alpha+1)), \mathbb{P}(z))} \mu(\b) \right)
  \le 
  \sum_{1 \le  \N(\alpha) \le u } \left(\sum_{\b \mid ((\alpha(t\alpha+1)), \mathbb{P}(z))} \lambda_{\b} \right)^2.
$$
Rearranging the sums we get
$$
  \sum_{1 \le \N(\alpha) \le u } \left(\sum_{\b \mid ((\alpha(t\alpha+1)), \mathbb{P}(z))} \lambda_{\b} \right)^2
  = 
  \sum_{\b_1, \b_2 \mid \mathbb{P}(z) \atop{\N\b_i \le z}} \lambda_{\b_1} \lambda_{\b_2}   \sum_{1 \le \N(\alpha) \le u \atop{ \b_i \mid (\alpha(t\alpha+1)) }} 1.
$$
Note that $\b = [\b_1, \b_2]$ is squarefree. If $\b_{(t)}$ is the maximal divisor of $\b$ co-prime to $(t)$,
we need to count $\alpha$ in $2^{\Omega_{\K}(\b_{(t)})}$ classes in $\O_{\K}/\b$.
Using \thmref{counting} for $\a = \O_{\K}$, $\q = \b$, we get
for $ z \le u^{\frac{1}{2}}$
$$
 \sum_{\b_1, \b_2 \mid \mathbb{P}(z) \atop{\N\b_i \le z}} \lambda_{\b_1} \lambda_{\b_2} \sum_{1 \le \N(\alpha) \le u \atop{ [\b_1, \b_2] \mid (\alpha(t\alpha+1)) }} 1
 =
  \sum_{\b_1, \b_2 \mid \mathbb{P}(z) \atop{\N\b_i \le z}} \lambda_{\b_1} \lambda_{\b_2} \left(\frac{c_\K  2^{\Omega_{\K}([\b_1,\b_2]_{(t)})} u}{\N[\b_1,\b_2]} + \rO^*\left( 10^{14} 2^{\Omega_\K([\b_1,\b_2]_{(t)})} \sqrt{\frac{u}{\N[\b_1,\b_2]}} \right) \right),
$$
where $c_{\K}$ is the constant given by $\frac{2\pi}{\sqrt{|d_{\K}|}} $.
We first deal with the error term.
By \lemref{errorterm}, we get
$$
  \sum_{\b_1, \b_2 \mid \mathbb{P}(z) \atop{\N\b_i \le z}} |\lambda_{\b_1} \lambda_{\b_2}| \frac{2^{\Omega_\K([\b_1,\b_2]_{(t)})}}{\sqrt{\N[\b_1,\b_2]}}
  \le 
  3^{32} \cdot z.
$$
 The main term now becomes $
\displaystyle 
 \sum_{\b_1, \b_2 \mid \mathbb{P}(z)} \frac{\lambda_{\b_1} \lambda_{\b_2} }{f_{(t)}([\b_1,\b_2])}$
 which in turn is
$$
\sum_{\b_1, \b_2 \mid \mathbb{P}(z)} \frac{\lambda_{\b_1} \lambda_{\b_2} f_{(t)}((\b_1,\b_2)) }{f_{(t)}(\b_1)f_{(t)}(\b_2)}
=
\sum_{\b_1, \b_2 \mid \mathbb{P}(z)} \frac{\lambda_{\b_1} \lambda_{\b_2}}{f_{(t)}(\b_1)f_{(t)}(\b_2)} \sum_{\a \mid (\b_1,\b_2)} f_{1,(t)}(\a)
=
\sum_{\a \mid \mathbb{P}(z)} f_{1,(t)}(\a) \left(\sum_{ \c \mid \mathbb{P}(z) \atop{\a \mid \c }} \frac{\lambda_{\c}}{f_{(t)}(\c)}\right)^2.
$$
Further, we observe that
\begin{equation*}
\sum_{ \c \mid \mathbb{P}(z) \atop{\a \mid \c }} \frac{\lambda_{\c}}{f_{(t)}(\c)} 
~=~  
S_{\O_\K}(z)^{-1} \sum_{\c \mid \mathbb{P}(z) \atop{\a \mid \c}}\frac{\mu(\c)}{f_{1,(t)}(\c)} 
\sum_{1 \le \N(\g) \le \frac{z}{\N(\c)} \atop{(\g, (2310)\c) = \O_{\K}}} \frac{\mu^2(\g)}{f_{1,(t)}(\g)} .
\end{equation*}
Writing $\c = \h\a$ with $(\h,\a) = \O_{\K}$, we get
\begin{eqnarray*}
&& \frac{\mu(\a)}{f_{1,(t)}(\a)} S_{\O_\K}(z)^{-1}   \sum_{1 \le\N(\h) \le \frac{z}{\N(\a)} \atop{\h \mid \mathbb{P}(z) \atop{(\h,\a) = \O_{\K}} } } 
\frac{\mu(\h)}{f_{1,(t)}(\h)} \sum_{1 \le \N(\g) \le \frac{z}{\N(\h\a)} \atop{(\g,  (2310)\h\a) = \O_{\K}}} 
\frac{\mu^2(\g)}{f_{1,(t)}(\g)} \\
& = & \frac{\mu(\a)}{f_{1,(t)}(\a)} S_{\O_\K}(z)^{-1} \sum_{1 \le \N(\h) \le \frac{z}{\N(\a)} \atop{\h \mid \mathbb{P}(z) \atop{(\h,\a) = \O_{\K}}} }  
\sum_{1 \le \N(\g) \le \frac{z}{\N(\h\a)} \atop{(\g, (2310)\h\a) = \O_{\K}}} \mu(\h)\frac{\mu^2(\g\h)}{f_{1,(t)}(\g\h)}.\\
\end{eqnarray*}
Setting $\a_1 = \g\h$ gives
\begin{equation*}
\sum_{\c \mid \mathbb{P}(z) \atop{\a \mid \c}}\frac{\lambda_{\c}}{f_{(t)}(\c)}
~=~ 
\frac{\mu(\a)}{f_{1,(t)}(\a)} S_{\O_\K}(z)^{-1} \sum_{1 \le \N(\a_1) \le \frac{z}{\N(\a)} \atop{(\a_1, (2310)\a) = \O_{\K} \atop{\a_1 \mid \mathbb{P}(z)}}}
\frac{\mu^2(\a_1)}{f_{1,(t)}(\a_1)} \sum_{\h\mid \a_1} \mu(\h)
~ = ~ 
S_{\O_\K}(z)^{-1}  \frac{\mu(\a)}{f_{1,(t)}(\a)}.
\end{equation*}
Combining the above, the main term is
$
c_{\K} u   S_{\O_{\K}}(z)^{-1}.
$
However,
\begin{eqnarray*}
 S_{\O_{\K}}(z)
& = &
\sum_{1 \le \N(\a) \le z \atop{(\a, (2310)) = \O_{\K}}} \frac{\mu^2(\a)}{f_{1,(t)}(\a)}
=
\sum_{1 \le \N(\a) \le z \atop{(\a, (2310)) = \O_{\K}}} \frac{\mu^2(\a)}{f_{(t)}(\a)} \cdot \frac{f_{(t)}(\a)}{f_{1,(t)}(\a)} \\
& = & 
\sum_{1 \le \N(\a) \le z \atop{(\a, (2310)) = \O_{\K}}} \frac{\mu^2(\a)}{f_{(t)}(\a)} \cdot \prod_{\P \mid \a} \frac{f_{(t)}(\P)}{ \mu(\P) + f_{(t)}(\P)} \\
& \ge &
\sum_{1 \le \N(\a) \le z \atop{(\a, (2310)) = \O_{\K}}}  \frac{1}{f_{(t)}(\a)} 
-
\sum_{1 \le \N(\a) \le z \atop{ \a \text{ non-sq-free} \atop{(\a, (2310)) = \O_{\K}}}}  \frac{1}{f_{(t)}(\a)} .
\end{eqnarray*}
Consider
$$
\sum_{1 \le \N(\a) \le z \atop{ \a \text{ non-sq-free} \atop{(\a, (2310)) = \O_{\K}}}}  \frac{1}{f_{(t)}(\a)}
\le 
\sum_{1 \le \N\b \le z \atop{(\b, (2310)) = \O_{\K}}}  \sum_{\N\P \le \sqrt{z} \atop{ (\P , (2310)) = \O_{\K}}} \frac{1}{f_{(t)}(\b)f_{(t)}(\P)^2}
=
\sum_{1 \le \N\b \le z \atop{ (\b , (2310)) = \O_{\K}}} \frac{1}{f_{(t)}(\b)} \left( \sum_{\N\P \le \sqrt{z} \atop{ (\P , (2310)) = \O_{\K}}} \frac{1}{f_{(t)}(\P)^2} \right).
$$
Looking at the sum
$$
 \sum_{\N\P \le \sqrt{z} \atop{ (\P , (2310)) = \O_{\K}}} \frac{1}{f_{(t)}(\P)^2} \le   \sum_{\N\P \le \sqrt{z} \atop{ (\P , (2310)) = \O_{\K}}} \frac{4}{\N(\P)^2} 
 \le
 \frac{8}{11}.
$$

Therefore,
\begin{equation}\label{Sok}
S_{\O_\K}(z) \ge \frac{3}{11} \sum_{1 \le \N(\a) \le z \atop{(\a, (2310)) = \O_{\K}}}  \frac{1}{f_{(t)}(\a)} .
\end{equation}
Let us consider
$$
\sum_{1 \le \N\a \le z \atop{(\a, (2310)) = \O_{\K}}}  \frac{1}{f_{(t)}(\a)} 
= 
\sum_{1 \le \N\a \le z \atop{(\a, (2310)) = \O_{\K}}}  \frac{2^{\Omega_{\K}(\a_{(t)})}}{\N\a}
\ge
\sum_{1 \le \N\a \le z \atop{(\a, (2310)) = 
\O_{\K}}}  \frac{d_{(t)}(\a)}{\N\a}.
$$
Here  $d_{(t)}(\a)$ is the number of ideals dividing $\a$ and co-prime to $(t)$.
We now look at
\begin{eqnarray*}
\prod_{\P \mid (t)} \left(1 - \frac{1}{\N\P}\right)^{-1} \sum_{1 \le \N\a \le z \atop{(\a, (2310)) = \O_{\K}}}  \frac{d_{(t)}(\a)}{\N\a}
& = &
\sum_{1 \le \N\a \le z \atop{(\a, (2310)) = \O_{\K}}}  \frac{d_{(t)}(\a)}{\N\a} \sum_{\P \mid \b \Rightarrow \P \mid (t)} \frac{1}{\N\b} \\
& = &
\sum_{(0) \neq \b \subset \O_{\K}} \frac{1}{\N\b} \sum_{1 \le \N\a \le z, \a \mid \b \atop{(\a, (2310)) = \O_{\K} \atop{\P \mid \frac{\b}{\a} \Rightarrow \P \mid (t)}}}  d_{(t)}(\a) \\
& \ge &
\sum_{\b \subset \O_{\K} \atop{1 \le \N\b \le z \atop{(\b, (2310)) = \O_{\K}}}} \frac{d(\b)}{\N\b}. \\
\end{eqnarray*}
In order to compute the above sum, we consider the following 
\begin{eqnarray*}
\sum_{\b \subset \O_{\K} \atop {1 \le \N\b \le z \atop{(\b, (2310)) = \O_{\K}}}} d(\b)
& = &
\sum_{\b \subset \O_{\K} \atop{1 \le \N\b \le z \atop{(\b, (2310)) = \O_{\K}}}} \sum_{\e \mid \b} 1
=
\sum_{\N\e \le z \atop{(\e, (2310))=1}} \sum_{ \b \subset \O_{\K} \atop{1 \le \N\b \le z \atop{\e \mid \b, (\b, (2310)) = \O_{\K}}}}1
=
\sum_{\N\e \le z \atop{(\e,(2310))=1}} \sum_{ \b \subset \O_{\K} \atop{1 \le \N\b \le  \frac{z}{\N\e} \atop{ (\b, (2310)) = \O_{\K}}}}1.
\end{eqnarray*}
We have by \thmref{asymfinal}
$$
 \sum_{ \b \subset \O_{\K} \atop{1 \le \N\b \le  \frac{z}{\N\e} \atop{(\b, (2310)) = \O_{\K}}}}1
 =
 \frac{\rho_{\K} \varphi((2310))z}{2310^2\N\e}
 +
 O^*\left( 3 \cdot 10^{24} h_{\K}  \left(\frac{z}{\N\e}\right)^{\frac{1}{2}}
 +
  3 \cdot 10^{12}  h_{\K}  \right).
$$
Therefore, we have
$$
 \sum_{ \b \subset \O_{\K}
  \atop{1 \le \N\b \le  \frac{z}{\N\e}
  \atop{(\b, (2310)) = \O_{\K} } }}1
\ge
 \frac{\rho_{\K} \varphi((2310))z}{2310^2\N\e}
 -
 3 \cdot 10^{24} h_{\K} \left(\frac{z}{\N\e}\right)^{\frac{1}{2}}
 -
 3 \cdot 10^{12} h_{\K}.
$$
By partial summation formula,
\begin{eqnarray*}
\sum_{\e \subset \O_{\K} \atop{ \N\e \le z \atop{(\e,(2310))=\O_{\K}}}} \frac{1}{\N\e}
\ge
\left( \frac{\rho_{\K} \varphi((2310))}{2310^2} \right) \log ez
- 7 \cdot 10^{24} h_{\K}.
\end{eqnarray*}
Similarly, again by partial summation, we get
\begin{eqnarray*}
\sum_{\e \subset \O_{\K} \atop{ \N\e \le z \atop{(\e,(2310))=\O_{\K}}}} \frac{1}{\sqrt{\N\e}}
& \le &
\frac{2\rho_{\K} \varphi((2310))}{2310^2}  \sqrt{z} + 6 \cdot 10^{24} h_{\K} \log z.
\end{eqnarray*}
Combining the above inequalities,
\begin{eqnarray*}
\sum_{\b \subset \O_{\K} \atop {1 \le \N\b \le z \atop{(\b, (2310)) = \O_{\K}}}} d(\b)
& \ge &
\sum_{\N\e \le z \atop{(\e,(2310))=\O_\K}} 
\left(  \frac{\rho_{\K} \varphi((2310))z}{2310^2\N\e}
 -
 4 \cdot 10^{24} h_{\K} \left(\frac{z}{\N\e}\right)^{\frac{1}{2}} \right)\\
& \ge &
\left(  \frac{\rho_{\K} \varphi((2310))}{2310^2}
 \right)^2 z \log ez
 -
 4.81 \cdot 10^{49} \rho_{\K} h_{\K} \sqrt{|d_{\K}|} z
 \\
  & \ge &  
\frac12 \left(  \frac{\rho_{\K} \varphi((2310))}{2310^2}
 \right)^2 z \log ez,
\end{eqnarray*}
for 
$\log ez \ge \frac{ 2310^2 \cdot 10^{50}  h_{\K}  \sqrt{|d_{\K}|}}{\rho_{\K} \varphi((2310))^2}$.
Now, using the fact that $h_{\K} \ge \rho_{\K} \sqrt{|d_{\K}|}/ (2\pi)$ and that $\varphi((2310))~\ge~480^2$,
we get
$\log ez \ge 2 \cdot 10^{45} |d_{\K}|$.
We now apply partial summation once again to get
$$
\sum_{\b \subset \O_{\K} \atop{1 \le \N\b \le z\atop{(\b, (2310)) = \O_{\K}}}} \frac{d(\b)}{\N\b}
\ge 
  \frac{\left(6\rho_{\K}\right)^2}{2 \cdot 10^{10}}
  \log^2 z ~~~
 \text{ for }~~~
 \log ez \ge 2 \cdot 10^{45} |d_{\K}|.
$$
Combining the above for
$  \log ez \ge 2 \cdot 10^{45} |d_{\K}|$,
$$
S_{\O_{\K}}(z)^{-1} \le \prod_{\P \mid (t)} \left(1 - \frac{1}{\N\P}\right)^{-1} \frac{22 \cdot 10^{10}}{ 3 \left(  6 \rho_{\K} 
 \right)^2 \log^2 z }.
$$
If we suppose further that $z \le \sqrt{u}$, we get
\begin{eqnarray*}
\sum_{1 \le \N(\alpha) \le u \atop{\alpha, t\alpha+1 ~\in~ \p_{\K}}} 1 
& \le & 4 |\mu_{\K}|z 
+  
\frac{2\pi}{\sqrt{|d_{\K}|}} \prod_{\P \mid (t)} \left(1 - \frac{1}{\N\P}\right)^{-1} \frac{22 \cdot 10^{10} u}{ 3\left( 6 \rho_{\K}
 \right)^2 \log^2 z } 
 + 
10^{14} \cdot 3^{32} \cdot  z  \sqrt{u}\\
 & \le & 
\frac{2\pi}{\sqrt{|d_{\K}|}} \cdot  \frac{8 \cdot 10^{10} \N(t) u}{\varphi((t)) \left(  6\rho_{\K}
 \right)^2 \log^2 z } 
 + 
 10^{30}z\sqrt{u}.\\
\end{eqnarray*}
Further, since we have $\log z \le 2\sqrt{z}$, we get for $ \displaystyle z = \frac{u^{1/4}}{2}$,
\begin{eqnarray*}
\sum_{1 \le \N(\alpha) \le u \atop{\alpha, t\alpha+1~\in~ \p_{\K}}} 1 
 & \le &
\left(  \frac{2 \cdot 10^{10}}{ \rho_{\K}^2 \sqrt{|d_{\K}|}    } 
  \right) \frac{\N(t) u}{\varphi((t)) \log^2 \frac{ u^{1/4}}{2}}
  +
  10^{30} u^{3/4}.
\end{eqnarray*}
\end{proof}

\begin{lem}\label{comparisonterm}
For $Q \ge \exp(10^{45}|d_{\K}|)$, we have
$$
h_{\K} \pi^*(Q) \ge \frac{2Q}{25\log Q}.
$$
\end{lem}
\begin{proof}
	We first observe that for $ y \ge 8$, we have 
	$$\int_2^y \frac{1}{\log t} dt  \ge \frac{y}{\log y}.
	$$
	For $y \ge 200^5$, we have $4/\log 2 \le y^{1/5}/\log y$ and this implies
	$$
	\int_2^y \frac{1}{\log t} dt \le \int_2^{y^{4/5}} \frac{1}{\log 2} dt +  \int_{y^{4/5}}^y \frac{1}{\log y^{4/5}} dt
	\le
	\frac{y^{4/5}}{\log 2} + \frac{5y}{4\log y}
	\le
	\frac{3y}{2\log y}.
	$$
We have
$$
\pi^*(Q) = \pi(Q, \K(\O_\K)/\K, \sigma_0) - \pi \left(Q/2, \K(\O_\K)/\K, \sigma_0\right),
$$
where $\sigma_0$ is the trivial element
of the Galois group of $\K(\O_\K)/\K$.
We now get by \corref{principalprimecounting},
\begin{eqnarray*}
h_{\K}\pi^*(Q) & \ge & \frac{Q}{\log Q}- 5 h_\K \sqrt{Q} \log |d_{\K}| - 2 h_\K\sqrt{Q} \left(\frac{\log Q}{8\pi} + 9 \right) \\
&-&  \frac{3Q}{4\log Q/2} - 5  h_\K \sqrt{\frac{Q}{2}} \log |d_{\K}| - 2  h_\K \sqrt{\frac{Q}{2}} \left(\frac{\log \frac{Q}{2}}{8\pi} + 9 \right).
\end{eqnarray*}
Further since $Q \ge 2^{6}$, we get
\begin{eqnarray*}
h_{\K}\pi^*(Q) & \ge & \frac{Q}{\log Q}- 5 h_\K \sqrt{Q} \log |d_{\K}| - 2 h_\K\sqrt{Q} \left(\frac{\log Q}{8\pi} + 9 \right) \\
&-&  \frac{18Q}{20\log Q} - 5  h_\K \sqrt{\frac{Q}{2}} \log |d_{\K}| - 2  h_\K \sqrt{\frac{Q}{2}} \left(\frac{\log \frac{Q}{2}}{8\pi} + 9 \right).
\end{eqnarray*}
Simplifying the above and using \lemref{DGRS} and \lemref{size}, we get
\begin{eqnarray*}
h_{\K}\pi^*(Q) 
& \ge &
\frac{Q}{\log Q} \left(\frac{1}{10} - \frac{ 10^3  |d_{\K}|^{3/4} \log |d_{\K}| \log Q}{\sqrt{Q}} - \frac{400  |d_{\K}|^{3/4} \log Q}{\sqrt{Q}} \left(\frac{\log Q}{8\pi} + 9 \right)  \right).
\end{eqnarray*}
We now bound each of the negative terms separately using the lower bound for $Q$ as follows.
Using the fact that $\log x/\sqrt{x}$ is a decreasing function for $x \ge 8$, for the first term, we get
$$
 \frac{  10^3  |d_{\K}|^{3/4} \log |d_{\K}| \log Q}{\sqrt{Q}} \le \frac{ 10^{48} |d_{\K}|^{7/4} \log |d_{\K}|}{ \exp(\frac{10^{45}|d_{\K}|}{2})} 
 \le \frac{1}{100}.
$$
Similarly, for the second term,
\begin{eqnarray*}
\frac{400  |d_{\K}|^{3/4}\log Q}{\sqrt{Q}} \left(\frac{\log Q}{8\pi} + 9 \right)
~ \le~ 
\frac{400 \cdot 10^{90} |d_{\K}|^{11/4}}{8\pi \exp(\frac{10^{45} |d_{\K}|}{2})} + \frac{3600 \cdot 10^{45} |d_{\K}|^{7/4}}{ \exp(\frac{10^{45} |d_{\K}|}{2})} ~\le~ \frac{1}{100}.\\
\end{eqnarray*}
This gives us the lemma.
\end{proof}

Next two lemmas are devoted to the proof of \propref{zeroes}. We first compute an upper bound of the real part of the function $\frac{\Gamma'_\chi}{\Gamma_\chi}(s)$ for $s \in \C$ (see equation \eqref{eq100} for its definition).

\begin{lem}\label{regamma}
		For $ \Re(s)> 1$ and number field $\K$, we have
		$$
		\Re\left( \frac{\Gamma'_\chi}{\Gamma_\chi}(s)\right) < 
		n_\K  \left(\log\left(\frac{|s+1|}{2} + 2\right) - \frac{\log \pi}{2} \right).
		$$
	\end{lem}
	
	\begin{proof}
		Taking the logarithmic derivative on both sides of \eqref{eq100}, we obtain
		$$
		\frac{\Gamma'_\chi}{\Gamma_\chi}(s) = - \frac{n_\K}{2} \log \pi + \frac{a_{\chi}}{2} \frac{\Gamma'}{\Gamma}\left( \frac{s+1}{2}\right) + \frac{(n_{\K} - a_{\chi})}{2} \cdot \frac{\Gamma'}{\Gamma}\left( \frac{s}{2}\right).
		$$
		Using \lemref{ahnlem}, for $\Re(s)>1$, we get
		\begin{align*}
		\Re \left( \frac{\Gamma'_\chi}{\Gamma_\chi}(s)\right) & < - \frac{n_\K}{2} \log \pi +  a_{\chi} \log\left(\frac{|s+1|}{2} + 2\right)  + (n_\K -a_{\chi}) \log\left(\frac{|s|}{2} + 2\right) \\
		 & < n_\K  \left(\log\left(\frac{|s+1|}{2} + 2\right) - \frac{\log \pi}{2} \right).
		\end{align*}
	\end{proof}

The following lemma gives us an upper bound of the absolute value of the logarithmic derivative of Hecke $L$- functions.
		\begin{lem}\label{lprime}
			Let $L(s,\chi)$ be the 
			Hecke $L$-function corresponding to a generalized
			Dirichlet character $\chi$ modulo integral ideal $\c$.
			Then for $ \Re(s)= \sigma > 1$ and number field $\K$, we have
			$$
			\Big|\frac{L'}{L}(s, \chi) \Big| < \frac{n_\K}{\sigma - 1}.
			$$
		\end{lem}

	\begin{proof}
		Taking logarithmic derivative of Euler product of 
		Hecke L function and Dedekind zeta function, we obtain
			$$
			\Big|\frac{L'}{L}(s, \chi) \Big| = \Big |\sum_{\P \nmid \c \atop{ \P \neq (0)}} \frac{ \log \N \P}{ \frac{(\N \P)^s}{ \chi([\P])} - 1}\Big|  \le \sum_{\P \neq (0)} \frac{ \log \N \P}{(\N \P)^{\sigma} - 1}
			= - \frac{\zeta'_\K}{\zeta_\K}(\sigma).
			$$
		Similarly, we note that 
		$$
		- \frac{\zeta'_\K}{\zeta_\K}(\sigma) \le - n_\K \frac{\zeta'}{\zeta}(\sigma) < \frac{n_\K}{\sigma - 1}, \phantom{M} \sigma > 1.
		$$
		The last inequality follows from \lemref{halllem}. Now
		combining the above identities, we get
		$$
		\Big|\frac{L'}{L}(s, \chi) \Big| < \frac{n_\K}{\sigma - 1}.
		$$
	\end{proof}
	
	Now let $L(s,\chi)$ be the 
	Hecke $L$-function corresponding to  generalized
	Dirichlet character $\chi$ modulo $\c$ and $\c'$ be the conductor of $\chi$.
Then following  Lagarias and Odlyzko \cite[Lemma 5.4]{LO}, we have
	the following result which gives us an explicit bound on the number of zeros
	of $L(s,\chi)$.
	
	\begin{prop}\label{zeroes}
		For a real number $t$, let $N_{t, \chi}$ denote the number of zeros 
	  (denoted 
		$\rho = \alpha + i\beta$)
		 of $L(s,\chi)$
		with $0 <\alpha < 1$ and $|\beta - t| \le 1$.
		Then
		$$
		N_{t, \chi} < 5 (3(n_\K + 1) + \log ( |d_\K| \N\c') +  2n_\K \log( |t| + 2)).
$$
Moreover, when $\K \neq \Q$ then
$$
N_{t, \chi} < 50 n_\K \log (|d_\K| \N\c'~(|t|+ 2)).
$$
\end{prop}
	
\begin{proof}
From \lemref{diskolem}, we have
\begin{align*}
\Re \left(\sum_\rho \left(\frac{1}{s - \rho} + \frac{1}{s- \overline{\rho}}\right)\right) &=  \Re\left(\frac{L'}{L}(s, \chi)\right) + \Re\left(\frac{L'}{L}(s, \overline{\chi})\right) +  \log ( |d_\K| \N\c') \\ 
&\, + \mathbf{1}_{\chi}\Re \left(\frac{2}{s} + \frac{2}{s-1}\right) + 2 \Re \left(\frac{\Gamma'_\chi}{\Gamma_\chi}(s)\right).
\end{align*}
Evaluating it
at $s=2+it$ and using \lemref{regamma} and \lemref{lprime}, we obtain
\begin{align}\label{pink}
\Re \left(\sum_\rho \left(\frac{1}{s - \rho} + \frac{1}{s- \overline{\rho}} \right)\right) & < n_\K (2 - \log \pi) +  \log ( |d_\K| \N\c') + 2 \Re \left(\frac{1}{2+it} + \frac{1}{1+it}\right) \nonumber \\
&  + 2 n_\K  \log\left(\frac{|3+it|}{2} + 2\right) \nonumber \\
& < 3 n_\K + 3 + \log ( |d_\K| \N\c') +  2n_\K \log( |t| + 2).
\end{align}
Clearly, when $\K \neq \Q$, the above quantity is bounded above by
$ 10 n_\K \log (|d_\K| \N\c'~(|t|+ 2)).$\\
As noticed by Lagarias and Odlyzko \cite[Lemma 5.4]{LO}, we have
\begin{align}\label{pinkk}
\Re \left(\sum_\rho \left(\frac{1}{s - \rho} + \frac{1}{s- \overline{\rho}} \right)\right)  \geq \sum_{\rho \atop{|\beta - t| \leq 1}} \frac{2-\alpha}{(2-\alpha)^2 + (t-\beta)^2} > \frac{1}{5} N_{t, \chi},
\end{align}
since $\alpha$ is between 0 and 1. Substituting \eqref{pinkk} into \eqref{pink}, we get our desired results.
\end{proof}

\section{Bounds on the average value of the Euler-Kronecker constant}

Throughout this section, we will use the following notations. Let $\K$ be an imaginary quadratic number field and $\q$ be a non-zero principal prime ideal of $\K$. Let $\chi$ be a generalized Dirichlet character $\mod \q$. 
Now after summing over all non-trivial characters in the following function (see equation \eqref{eq101})
$$\Phi_\chi(x)=\frac{1}{x-1}\int_1^x \Big(\sum_{\a \subset \mathcal{O}_\K \atop 1 \leq \N\mathfrak{a} \leq t}  \frac{\Lambda(\mathfrak{a})}{\N\mathfrak{a}}\chi([\mathfrak{a}])\Big)dt,$$
 we get by orthogonality of characters of a finite group,
\begin{equation}\label{Phi}
	\sum_{\chi \neq \chi_0}\Phi_\chi(x)=\frac{1}{x-1}\int_1^x\Big(|H_{\mathfrak{q}}(\K)|
	\sum_{\a \subset \O_\K \atop {1 \leq \N\mathfrak{a} \leq t  \atop \mathfrak{a} \in [1]_{\q}} }
	\frac{\Lambda(\mathfrak{a})}{\N\mathfrak{a}}
	-\sum_{\a \subset \O_\K \atop {1 \leq \N\mathfrak{a} \leq t \atop (\mathfrak{a},\mathfrak{q})=1}}\frac{\Lambda(\mathfrak{a})}{\N\mathfrak{a}}\Big)dt,
\end{equation}
where $[1]_{\q}$ is the trivial element
of the group $H_{\q}(\K)$.
Let 
\begin{eqnarray*} \psi(x,\q ,[1]_{\q})=\sum_{\a \subset \O_\K \atop  {1 \leq \N\mathfrak{a} \leq x\atop \mathfrak{a} \in [1]_{\q}}}\Lambda(\mathfrak{a})
&\textrm{and}&
\psi(x)=\sum_{\a \subset \O_\K \atop {1 \leq \N\mathfrak{a} \leq x}}\Lambda(\mathfrak{a}).
\end{eqnarray*}
Since the Galois group of $\K(\q)/\K$
is isomorphic to the group $H_{\q}(\K)$
under the Artin map, we can identify
the trivial element of the Galois group
of $\K(\q)/\K$ with the trivial element
$ [1]_{\q} \in H_{\q}(\K)$.
Therefore,
$\psi(x,\q,[1]_{\q})=\psi(x,\K(\q)/\K, \sigma_o)$
where $\sigma_o$ is the identity 
element of the Galois group of $\K(\q)/\K$. 
For a non-zero prime ideal $\mathfrak{q}$,
\begin{eqnarray}\label{coprimesum}
\sum_{\a \subset \O_\K \atop { 1 \leq \N\mathfrak{a} \leq x \atop (\mathfrak{a},\mathfrak{q})=1}}\frac{\Lambda(\mathfrak{a})}{\N\mathfrak{a}} \ge \sum_{\a \subset \O_\K \atop 1 \leq \N\mathfrak{a} \leq x}\frac{\Lambda(\mathfrak{a})}{\N\mathfrak{a}}- \frac{\log\N\mathfrak{q}}{\N\mathfrak{q}-1}.
\end{eqnarray}
By Abel's summation formula, we get
\begin{eqnarray*} 
\sum_{\a \subset \O_\K \atop {1 \leq \N\mathfrak{a} \leq x \atop \mathfrak{a} \in [1]_{\q}}}\frac{\Lambda(\mathfrak{a})}{\N\mathfrak{a}}=\frac{\psi(x,\q,[1]_{\q})}{x}+\int_1^x\frac{\psi(u,\q,[1]_{\q})}{u^2}du ~~\phantom{m}
\text{ and }~~ 
\sum_{\a \subset \O_\K \atop 1 \leq \N\mathfrak{a} \leq x}\frac{\Lambda(\mathfrak{a})}{\N\mathfrak{a}}= \frac{\psi(x)}{x}+\int_1^x\frac{\psi(u)}{u^2}du.
\end{eqnarray*}
Integrating over $t$ on both sides of \eqref{coprimesum}, we 
get
\begin{eqnarray*}
\int_1^x\sum_{\a \subset \O_\K \atop {1 \leq \N\mathfrak{a} \leq t \atop (\mathfrak{a},\mathfrak{q})=1}}\frac{\Lambda(\mathfrak{a})}{\N\mathfrak{a}}dt
&\geq& \int_1^x \frac{\psi(t)}{t}dt +\int_1^x\int_1^t\frac{\psi(u)}{u^2}du~dt
- \frac{\log\N\mathfrak{q}}{\N\mathfrak{q}-1}(x-1)\\
&\geq& \int_1^x \frac{\psi(t)}{t}dt +\int_1^x\int_u^x\frac{\psi(u)}{u^2}dt~du
- \frac{\log\N\mathfrak{q}}{\N\mathfrak{q}-1}(x-1)\\
&\geq&x \int_1^x \frac{\psi(u)}{u^2}du- \frac{\log\N\mathfrak{q}}{\N\mathfrak{q}-1}(x-1).
\end{eqnarray*}
Similarly, we
 get 
$$\int_1^x \sum_{\a \subset \O_\K \atop {1 \leq \N\mathfrak{a} \leq t \atop \mathfrak{a} \in [1]_{\q}}}\frac{\Lambda(\mathfrak{a})}{\N\mathfrak{a}}dt
=
x\int_1^x\frac{\psi(u,\mathfrak{q},[1]_{\q})}{u^2}du. $$ 
Hence, \eqref{Phi} 
becomes
\begin{equation}\label{Phi3}
	\sum_{\chi \neq \chi_0}\Phi_\chi(x)\leq\frac{x}{x-1}\int_1^x\left(|H_{\mathfrak{q}}(\K)|\frac{\psi(u,\mathfrak{q},[1]_{\q})}{u^2} - \frac{\psi(u)}{u^2}\right)du+ \frac{\log\N\mathfrak{q}}{\N\mathfrak{q}-1}.
\end{equation}

\begin{lem}\label{qbound}
We have
\begin{equation*}
\sum_{\qq}\frac{\log\N\mathfrak{q}}{\N\mathfrak{q}-1} \leq 14 + 4\frac{e^{75}|d_\K|^{1/3}(\log|d_\K|)^2}{\rho_\K},
	\end{equation*}
	where the sum runs over prime ideals $\q$ with $\N\q$ in the interval $(\frac{1}{2}Q, Q]$.
\end{lem}
\begin{proof}
Using \thmref{GLee}, we get
$$\left|\sum_{1 \le \N\q \leq x}\frac{\log\N\mathfrak{q}}{\N\mathfrak{q}} -\log x\right|\leq 3+ \frac{e^{75}|d_\K|^{1/3}(\log|d_\K|)^2}{\rho_\K} . $$
It implies that 
$$\sum_{\qq}\frac{\log\N\mathfrak{q}}{\N\mathfrak{q}-1} \leq \sum_{\qq}\frac{2\log\N\mathfrak{q}}{\N\mathfrak{q}}
\leq 2\log 2 + 12 + 4\frac{e^{75}|d_\K|^{1/3}(\log|d_\K|)^2}{\rho_\K} .
$$
\end{proof}
\begin{rmk}
We note, from the above, that we have 
$$
\sum_{\frac{Q}{2} \le \N\q \leq Q}\frac{\log\N\mathfrak{q}}{\N\mathfrak{q}-1}  < 8 \cdot 10^{45} |d_{\K}|.
$$
\end{rmk}

\begin{thm}\label{Phi1}
	Assuming GRH, for $x >2$ and $Q \geq 8 \exp(8 \cdot 10^{45} |d_{\K}|)$ , we have
	\begin{align*}
		\sideset{}{'}\sum_{\qq} \Big\lvert \sum_{\chi \neq \chi_0}\Phi_\chi(x)\Big\rvert  
&
\leq (2^{13} \cdot 10^{13} h_{\K} + 6000h_{\K}^2 + 10) \pi^*(Q)\log Q.
	\end{align*}
where the outer sum is over principal prime ideals and $ \pi^*(Q)$ denotes the number of principal prime ideals with norm in the interval $(\frac{1}{2}Q, Q]$.
\end{thm}

\begin{proof}
From now onwards, we will drop the prime over summation and always assume the sum to be over  principal prime ideals in the given interval. \\
	From  equation \eqref{Phi3}, for $x \geq 2$, we get
	\begin{equation}\label{1.2}
		\sum_{\qq} \Big\lvert \sum_{\chi \neq \chi_0}\Phi_\chi(x)\Big\rvert \leq \displaystyle 2\int_1^x\frac{\sum_{\qq}\Big\lvert|H_{\q}(\K)|\psi(u,\q,[1]_{\q})- \psi(u)\Big\rvert}{u^2}du + C_1(\K),
	\end{equation}
	where $\displaystyle C_1(\K)= 14 + 4\frac{e^{75}|d_\K|^{1/3}(\log|d_\K|)^2}{\rho_\K} .$\\
	For any principal ideal $\a \subset \O_{\K}$, let $\mathcal{T}_{\a}$ be the set of all principal prime ideals $\q$ such that there exists a generator $\pi_{\a}$ of $\a$ such that $\pi_{\a}\equiv 1 \mod \q$. Then, 
	\begin{equation}\label{1.1}
		\sum_{\qq}|H_{\q}(\K)|\psi(u,\q,[1]_{\q})=\sum_{\a \subset \O_\K \atop{1 \leq  \N(\mathfrak{a})\leq u \atop \a \text{ principal }}}\Lambda(\mathfrak{a})\sum_{\substack{\q \in \mathcal{T}_{\a} \atop \qq}}|H_{\q}(\K)|.
	\end{equation}
First we shall consider the case $u \leq Q^{5/2}$.\\
\textbf{Contribution from prime powers $\P^k$ for $k \geq 2$\;:}
 Consider a principal ideal $\a \subset \O_{\K}$, $1 \le \N\a \le u$ and let  $\mathfrak{p}_1, \ldots , \mathfrak{p}_j$ be the distinct prime ideals with norm in $(Q/2,Q] $ such that $\a \in [1]_{{\P}_i} $ for all $1 \leq i \leq j$. Therefore for each $1 \leq i \leq j$, we have $\pi_{\a,\P_i}$,  generators of $\a$, such that $\pi_{\a,\P_i} \equiv 1 \mod \P_i$. We observe that since $\a$ is a 
principal ideal, all the $\pi_{\a,\P_i}$
are associates. Then we have
$$
\prod_{i=1}^j \P_i \mid \prod_{w \in \mu_{\K}} (w\pi_{\mathfrak{a}, \P_1}-1).
$$
Let $c=|\mu_{\K}| $. 
Then 
\begin{eqnarray}\label{finitelymanyterms}
\frac{Q^j}{2^j}\leq \N(\mathfrak{p}_1 \ldots \mathfrak{p}_j) \leq \prod_{w \in \mu_{\K}} \N(w\pi_{\mathfrak{a}, \P_1}-1)\leq (4\N(\pi_{\mathfrak{a}, \P_1}))^{c}\leq (4u)^{c}\leq (4Q^{5/2})^{c},
\end{eqnarray}
where we have used the fact that $\N(\alpha - 1) \le 4\N(\alpha)$ for any non-zero integral element $\alpha$ in an imaginary quadratic field.
For $Q > 8$, comparing both sides of \eqref{finitelymanyterms}, we note that 
$$
2^{6j - 15c}< Q^{2j - 5c}\le 2^{2j + 4c}.
$$
This implies that $j \leq 5 |\mu_\K|$. So, the number of terms in the inner sum of \eqref{1.1} is bounded above by $5 |\mu_\K|$ for $Q > 8$.
Using \lemref{size}, we see that the contribution of powers of prime ideals to the sum \eqref{1.1} is 
	\begin{eqnarray*}
	\sum_{\a = \P^k \subset \O_\K \atop{1 \leq  \N\mathfrak{\P}^k \leq u \atop{ \P^k \text{ principal} \atop{k \ge 2} }}}\Lambda(\mathfrak{a})\sum_{\substack{\q \in \mathcal{T}_{\a} \atop \qq}}|H_{\q}(\K)|
	~ \leq~ (5|\mu_\K|)h_\K  Q\sum_{\substack{1 \leq \N\mathfrak{p}^k\leq u \atop k \geq 2}}\log\N\mathfrak{p}
	 ~\leq ~ 10 h_\K |\mu_\K|  Q\sqrt{u}\log u.
	\end{eqnarray*} 
Putting this in \eqref{1.2}, we get
\begin{eqnarray*}	
	2\int_1^{Q^{5/2}} \left( \sum_{\a = \P^k \subset \O_\K \atop{1 \leq  \N\mathfrak{\P}^k \leq u \atop{ \P^k \text{ principal} \atop{k \ge 2} }}}\Lambda(\mathfrak{a})\sum_{\substack{\q \in \mathcal{T}_{\a} \atop \qq}}|H_{\q}(\K)| \right) \frac{du}{u^2}
	 & \leq & 20 h_\K |\mu_\K|Q \int_1^{Q^{5/2}}\frac{\log u}{u^{3/2}} du    \leq  480 h_\K Q, 
\end{eqnarray*}
	where we use \eqref{CNF} and the fact that for $\Re(s)>0$,
\begin{equation}\label{INT}
\int_1^{\infty}\frac{\log^{n}t}{t^{s+1}}dt=\frac{\Gamma(n+1)}{s^{n+1}}.
\end{equation}	
Since $Q>\exp(10^{45} |d_{\K}|)$,
using \lemref{comparisonterm} and the fact that $\log |d_{\K}| \ge \log 2$, we get
\begin{eqnarray}\label{smalluprimepower}
2\int_1^{Q^{5/2}} \left( \sum_{\P^k \subset \O_\K \atop{1 \leq  \N\mathfrak{\P}^k \leq u \atop{ \P^k \text{ principal} \atop{k \ge 2} }}}\Lambda(\mathfrak{a})\sum_{\substack{\q \in \mathcal{T}_{\a} \atop \qq}}|H_{\q}(\K)| \right) \frac{du}{u^2}
~~\leq~~ 6000 h_\K^2 \pi^*(Q)\log Q.
\end{eqnarray}
\textbf{Contribution from primes $\P$\;:}
We now want to compute the contribution
of prime ideals to \eqref{1.1}.
To do so, we consider the sum
\begin{eqnarray*}
		\displaystyle
		\sum_{\P \subset \O_{\K} \atop{1 \leq \N\mathfrak{p}\leq u \atop \P \text{ principal }}}\log{\N\mathfrak{p} }\sum_{ \qq \atop \mathfrak{p} \in [1]_{\q}} |H_{\q}(\K)|.
\end{eqnarray*}
To bound this term for each $\P$ in the outer sum, we choose a generator $\pi_\P$
and count the number of ideals $\q$ such that $\pi_\P \equiv 1 \bmod \q$.
Therefore
$$
\displaystyle
		\sum_{\P \subset \O_{\K} \atop{1 \leq \N\mathfrak{p}\leq u \atop \P \text{ principal }}}\log{\N\mathfrak{p} }\sum_{\substack{\mathfrak{p} \in [1]_{\q} \atop \qq}} |H_{\q}(\K)|
		\le
\displaystyle
		\sum_{\pi_\P \in \O_{\K} \atop{1 \leq  \N(\pi_{\P})\leq u \atop \pi_\P  \in \p_\K}}\log{\N(\pi_{\P})}		\sum_{\substack{\pi_{\mathfrak{p}} \equiv 1 \bmod\q} \atop \qq} |H_{\q}(\K)|.
$$
For each principal prime ideal $\q$, we now
fix a generator $\pi_{\q}$. The condition
$\pi_{\mathfrak{p}} \equiv 1 \bmod\q$ is equivalent to the existence of a $\delta 
\in \O_{\K}$ with $\N(\delta) \le 8u/Q$ such that
$
\pi_{\mathfrak{p}} - 1 = \delta \pi_{\q}.
$
Therefore, we have
\begin{eqnarray}\label{1.4}
		\sum_{\P \subset \O_{\K} \atop{1 \leq \N\mathfrak{p}\leq u \atop \P \text{ principal }}}\log{\N\mathfrak{p} }\sum_{\substack{\mathfrak{p} \in [1]_{\q} \atop \qq}} |H_{\q}(\K)|
		& \le &
\sum_{\pi_\P \in \O_{\K} \atop{1 \leq  \N(\pi_{\P}) \leq u \atop \pi_\P \in \p_{\K}}} \sum_{\substack{\qq \atop{ \pi_{\mathfrak{p}} -1= \pi_{\q}\delta \atop{ \text{ for some }\delta\in \O_{\K} }}}} |H_{\q}(\K)|\log{\N \pi_\P} \nonumber\\
	& \le &
\sum_{\pi_\P \in \O_{\K} \atop{1 \leq  \N(\pi_{\P}) \leq u \atop \pi_\P \in \p_{\K}}} \sum_{\substack{\frac12 Q < \N(\pi_{\q}) \le Q \atop{ \pi_{\mathfrak{p}} -1= \pi_{\q}\delta \atop{ \text{ for some }\delta\in \O_{\K} }}}} |H_{(\pi_\q)}(\K)|\log{\N(\pi_\P)} \nonumber\\
& \le &  
		\sum_{\delta \in \O_{\K} \atop 1 \leq \N(\delta) \leq \frac{8u}{Q}}
		\sum_{\pi_{\P},\pi_{\q}\in \p_{\K} \atop{\N(\pi_{\P}) \leq u 
				\atop{\frac12 Q < \N(\pi_{\q}) \le Q \atop{ \pi_{\mathfrak{p}}-1=\delta\pi_{\q}}}}}
		|H_{(\pi_\q)}(\K)|\log{\N(\pi_\P)}.
\end{eqnarray}
We split the outer  summation as 
$$\sum_{\delta \in \O_{\K} \atop 1 \leq \N(\delta) \leq \frac{8u}{Q}} = \sum_{\delta \in \O_{\K} \atop 1 \leq \N(\delta) < \frac{u}{Q}} + \sum_{\delta \in \O_{\K} \atop \frac{u}{Q} \leq \N(\delta) \leq \frac{8u}{Q}}$$
Now we shall compute the following integrals 
 \begin{align}
 I_1=& \int_1^{Q^{5/2}} \displaystyle
\sum_{\delta \in \O_{\K} \atop 1 \leq \N(\delta) < \frac{u}{Q}}
\sum_{\substack{\pi_{\P},\pi_{\q} \in \p_\K \atop{\N(\pi_{\P}) \leq u 
\atop{\frac12 Q < \N(\pi_{\q}) \le Q \atop \pi_{\mathfrak{p}}-1=\delta\pi_{\q}}}}}|H_{(\pi_{\q})}(\K)|\log{\N(\pi_{\P})}\frac{du}{u^2}\label{I_1int}\\
 I_2=& \int_1^{Q^{5/2}} \displaystyle
\sum_{\delta \in \O_{\K} \atop \frac{u}{Q} \leq \N(\delta) \leq \frac{8u}{Q}}
\sum_{\substack{\pi_{\P},\pi_{\q} \in \p_\K \atop{\N(\pi_{\P}) \leq u 
\atop{\frac12 Q < \N(\pi_{\q}) \le Q \atop \pi_{\mathfrak{p}}-1=\delta\pi_{\q}}}}}|H_{(\pi_{\q})}(\K)|\log{\N(\pi_{\P})}\frac{du}{u^2}.\label{I_2int}
\end{align}
\emph{Estimating $I_1$ :}
For $\delta \in \O_{\K}\backslash\{0\}$, we note that counting the number of pairs of prime elements $(\pi_\P, \pi_\q)$ such that $\pi_\P-1=\delta\pi_\q$ in the inner sum in \eqref{I_1int} is same as counting number of $\alpha\in \p_\K$ such that $\delta\alpha+1$ also belongs to $\p_\K$.
Using the conditions $\pi_{\mathfrak{p}}-1=\delta\pi_{\q}$ and $Q/2~<~\N(\pi_\q)\leq~Q$, we get $\N(\pi_\P)\leq 4\N(\pi_\P -1) \leq 4\N(\delta)Q$.
Also, $\N(\alpha)=\N(\pi_\q) \le Q$.  Using  Proposition \ref{selbergsieve} we bound this term by

\begin{equation*}\label{sum2}
\sum_{\substack{\pi_{\P},\pi_{\q} \in \p_\K \atop{\N(\pi_{\P}) \leq u 
\atop{\frac12 Q < \N(\pi_{\q}) \le Q \atop \pi_{\mathfrak{p}}-1=\delta\pi_{\q}}}}} 1\le \sum_{1 \le \N(\alpha) \le Q \atop{\alpha, \delta\alpha+1 \in \p_\K}} 1 
  \le 
    \frac{2^5 \cdot 10^{10} \cdot \N(\delta) Q}{ \rho_\K^2 \sqrt{|d_{\K}|} \varphi((\delta)) \log^2 \frac{Q}{16}}
 +
 10^{30} Q^{3/4}.
\end{equation*}
Using equation \eqref{sum2} and \lemref{size}, we get

$$I_1 \leq h_\K Q \int_{1}^{Q^{5/2}} \sum_{\delta \in \O_{\K} \atop 1 \leq \N(\delta) < \frac{u}{Q}}\left( \frac{2^5 \cdot 10^{10} \cdot \N(\delta) Q}{ \rho_\K^2 \sqrt{|d_{\K}|} \varphi((\delta)) \log^2 \frac{Q}{16}}
 +
 10^{30} Q^{3/4}\right) \log 4Q\N(\delta) \frac{du}{u^2}. $$
	Interchanging sum and integral, we get
\begin{align*}
I_1 & \leq h_\K Q  \sum_{\delta \in \O_{\K} \atop 1 \leq \N(\delta) < Q^{3/2}}\left( \frac{2^5 \cdot 10^{10} \cdot \N(\delta) Q}{ \rho_\K^2 \sqrt{|d_{\K}|} \varphi((\delta)) \log^2 \frac{Q}{16}}
 +
 10^{30} Q^{3/4}\right) \log 4Q\N(\delta) \int_{Q\N(\delta)}^{\infty}\frac{du}{u^2} \\
 & \leq
 h_\K Q  \sum_{\delta \in \O_{\K} \atop 1 \leq \N(\delta) < Q^{3/2}}\left( \frac{2^5 \cdot 10^{10} \cdot \N(\delta) Q}{ \rho_\K^2 \sqrt{|d_{\K}|} \varphi((\delta)) \log^2 \frac{Q}{16}}
 +
 10^{30} Q^{3/4}\right) \frac{\log 4Q\N(\delta)} {Q\N(\delta)}.
\end{align*}
Let
$$
I_{11} =  h_\K Q  \sum_{\delta \in \O_{\K} \atop 1 \leq \N(\delta) < Q^{3/2}}\left( \frac{2^5 \cdot 10^{10} \cdot \N(\delta) Q}{ \rho_\K^2 \sqrt{|d_{\K}|} \varphi((\delta)) \log^2 \frac{Q}{16}} \right) 
 \frac{\log 4Q\N(\delta)} {Q\N(\delta)}
$$
and
$$
I_{12} = h_\K Q  \sum_{\delta \in \O_{\K} \atop 1 \leq \N(\delta) < Q^{3/2}}
10^{30} Q^{3/4}\frac{\log 4Q\N(\delta)} {Q\N(\delta)}.
$$
Computing the summation in the first term and using \lemref{lem12}, we get
\begin{eqnarray*}
I_{11} & \le & \frac{2^4 \cdot 5 \cdot 10^{10} h_{\K} Q\log (2Q)}{ \rho_\K^2 \sqrt{|d_{\K}|} \log^2 \frac{Q}{16}} \sum_{\delta \in \O_{\K} \atop 1 \leq \N(\delta) < Q^{3/2}} \frac{1}{ \varphi((\delta))} \\
& \le &
\frac{2^4 \cdot 5 \cdot 10^{10} h_{\K} Q\log (2Q)}{ \rho_\K^2 \sqrt{|d_{\K}|} \log^2 \frac{Q}{16}}
 \left(\zeta_\K(2)^2 \rho_\K \log 4Q^{3/2} +  \sqrt{3} \cdot 2 \cdot 10^{15}\zeta_{\K}(3/2)\zeta_\K(2) h_\K \right).
\end{eqnarray*}
For $Q >  2^{9}$
we have
$\log 2Q \le 2\log{\frac{Q}{16}}
$.
Further, since $h_{\K} \le \rho_{\K} \sqrt{|d_{\K}|}$, we get
\begin{eqnarray*}
I_{11} & \le &
2^8 \cdot 3 \cdot 10^{11} Q
 +  \frac{10^{40}h_{\K} Q }{\rho_{\K} \log \frac{Q}{16}}.
\end{eqnarray*}
Since $Q \ge 8\exp(8 \cdot 10^{45} \cdot |d_{\K}|)$, we have 
$\log Q/16 \ge 10^{45} \cdot |d_{\K}|$.
This gives us 
$$
I_{11}  \le 
2^8 \cdot 3 \cdot 10^{11} Q
 +  \frac{ Q }{\sqrt{|d_{\K}|}}.
$$
Using \lemref{comparisonterm}, we get
\begin{eqnarray}\label{i11}
I_{11} \le
\left( 2^5 \cdot 3 \cdot 10^{13} h_{\K} + \frac{1}{2} \right) \pi^*(Q) \log Q.
\end{eqnarray}
We now consider
$$
I_{12} \le \frac52 \cdot 10^{30} h_{\K} Q^{3/4} \log 2Q \sum_{\delta \in \O_{\K} \atop 1 \leq \N(\delta) < Q^{3/2}}
\frac{1} {\N(\delta)}.
$$
By partial summation formula, we have
\begin{equation}\label{sumndelta}
\sum_{\delta \in \O_{\K} \atop 1 \leq \N(\delta) < Q^{3/2}}
\frac{1} {\N(\delta)}
\le
2\rho_{\K} \log 2Q + 2 \cdot \sqrt{3} \cdot 10^{15} h_{\K}.
\end{equation}
We now substitute the bound \eqref{sumndelta} in bound for $I_{12}$.
Using \lemref{comparisonterm}, the
fact that $\log 2Q \le 16 (2Q)^{1/16}$,
and the bound $Q \ge 8\exp(8 \cdot 10^{45} \cdot |d_{\K}|)$
we get
\begin{equation}\label{i12}
I_{12} \le   10^{35}  \rho_\K h_{\K}^2    \frac{ \pi^*(Q) \log Q}{Q^{1/4}}
 +
 10^{49} h_{\K}^3 \frac{ \pi^*(Q) \log Q}{Q^{1/4}}
 \le
\frac12 \pi^*(Q) \log Q.
\end{equation}
Therefore,
\begin{equation}\label{i1}
I_1 \le (2^5 \cdot 3 \cdot 10^{13} h_{\K} +1)\pi^*(Q) \log Q.
\end{equation}
\emph{Estimating $I_2$ :}
In the same way as before, for $\delta \in \O_{\K}\backslash\{0\}$, counting the number of pairs of prime elements $(\pi_\P, \pi_\q)$ such that $\pi_\P-1=\delta\pi_\q$ in the inner sum in \eqref{I_2int} is same as counting number of $\alpha\in \p_\K$ such that $\delta\alpha+1$ also belongs to $\p_\K$. Using the condition $\pi_{\mathfrak{p}}-1=\delta\pi_{\q}$we get $\N(\delta\alpha)=\N(\delta\pi_\q)= \N(\pi_\P -1) \leq 4\N(\pi_\P)\leq 4u$.
 Using Proposition \ref{selbergsieve} we bound this term by 
\begin{equation}\label{sum}
\sum_{\substack{\pi_{\P},\pi_{\q} \in \p_\K \atop{\N(\pi_{\P}) \leq u 
\atop{\frac12 Q < \N(\pi_{\q}) \le Q \atop \pi_{\mathfrak{p}}-1=\delta\pi_{\q}}}}} 1\le \sum_{1 \le \N(\alpha) \le 4u/\N(\delta) \atop{\alpha, \delta\alpha+1 \in \p_\K}} 1 
\leq
 \frac{2^7 \cdot 10^{10} u}{ \rho_\K^2 \sqrt{|d_{\K}|} \varphi((\delta)) \log^2 \frac{u}{4\N(\delta)}}
 +
 10^{30} \left(\frac{4u}{ \N(\delta)} \right)^{3/4}.
\end{equation}
Using \lemref{size} and \eqref{sum}, we get  
\begin{align*}
	I_2 ~\leq~  \frac{2^7 \cdot 10^{10}h_\K  Q}{ \rho_{\K}^2 \sqrt{|d_{\K}|}{\log^2 (\frac{Q}{32})}}\int_1^{Q^{5/2}} \displaystyle
\sum_{\delta \in \O_{\K} \atop \frac{u}{Q} \leq \N(\delta) \leq \frac{8u}{Q}}\frac{u^{-1}\log u}{\varphi_\K((\delta))}du
~+~
3 \cdot 10^{30} h_{\K}Q 
\int_1^{Q^{5/2}} \displaystyle
\sum_{\delta \in \O_{\K} \atop \frac{u}{Q} \leq \N(\delta) \leq \frac{8u}{Q}}
\frac{ \log u } {\N(\delta)^{3/4}u^{5/4}} du.
\end{align*}	
We shall denote henceforth by
$$
I_{21} = \frac{2^7 \cdot 10^{10}h_\K  Q}{ \rho_{\K}^2 \sqrt{|d_{\K}|}{\log^2 ( \frac{Q}{32})}}\int_1^{Q^{5/2}} 
\sum_{\delta \in \O_{\K} \atop \frac{u}{Q} \leq \N(\delta) \leq \frac{8u}{Q}}\frac{u^{-1}\log u}{\varphi_\K((\delta))}du
$$
and
$$
I_{22} = 3 \cdot 10^{30} h_{\K}Q 
\int_1^{Q^{5/2}}
\sum_{\delta \in \O_{\K} \atop \frac{u}{Q} \leq \N(\delta) \leq \frac{8u}{Q}}
\frac{ \log u } {u^{5/4} \N(\delta)^{3/4}}du.
$$
Let us first consider $I_{21}$, 
interchanging sum and integral and using the bound $h_{\K} \le \rho_{\K} \sqrt{|d_{\K}|}$, we get

$$
I_{21} \le \frac{2^7 \cdot 10^{10} Q}{ \rho_{\K} \log^2 (\frac{Q}{32})}
\sum_{\delta \in \O_{\K} \atop 1 \leq \N(\delta) \leq 8Q^{3/2}}
\frac{1}{\varphi_\K((\delta))}
\int_{\frac{Q\N(\delta)}{8}}^{Q\N(\delta)}
\frac{\log u}{u}du.
$$
Integrating and using \lemref{lem12} and the fact that $\log (3h_{\K}) \le 8(3h_{\K})^{1/8}$, we get
\begin{eqnarray*}
I_{21} 
& \le &
\frac{2^6 \cdot 10^{10} \cdot \log 8 \cdot Q}{ \rho_{\K} \log^2 (\frac{Q}{32})}
\sum_{\delta \in \O_{\K} \atop 1 \leq \N(\delta) \leq 8Q^{3/2}}
\frac{1}{\varphi_\K((\delta))}
\log \left(\frac{Q^2 \N(\delta)^2}{8} \right)\\
& \le &
\frac{2^9 \cdot 10^{11} \cdot \log 8 \cdot Q \log^2 (2Q)}{ \log^2 (\frac{Q}{32})} 
+
\frac{2^{8} \cdot 3^4 \cdot 10^{26} \cdot Q \log (2Q)}{ \rho_{\K} \log^2 (\frac{Q}{32})}
(h_{\K}\log(3h_{\K}))^{1/2} \\
& \le &
\frac{2^9 \cdot 10^{11} \cdot \log 8 \cdot Q \log^2 (2Q)}{ \log^2 (\frac{Q}{32})} 
+
\frac{2^{11} \cdot 3^4 \cdot 10^{26} \cdot Q \log (2Q)}{ \rho_{\K} \log^2 (\frac{Q}{32})}
h_{\K}^{9/16}. \\
\end{eqnarray*}

Using the fact that $2\log (Q/32) \ge  \log (2Q)$ for $Q \ge 2^{11}$, we get
$$
I_{21} \le 
2^{13} \cdot 10^{11} \cdot \log 8 \cdot Q 
+
\frac{2^{12} \cdot 3^4 \cdot 10^{26} Q}{ \rho_{\K} \log (\frac{Q}{32})} h_{\K}^{9/16}.
$$
Since $Q \ge 8\exp(8 \cdot 10^{45} |d_{\K}|)$ and bounding $h_{\K}$, we get
$$
I_{21} \le 
2^{13} \cdot 10^{11} \cdot \log 8 \cdot Q 
+
\frac{Q}{2 \cdot 10^{10} \rho_{\K}^{7/16}  |d_{\K}|^{23/32}}.
$$
We now apply \lemref{comparisonterm}
and we bound $\rho_{\K}$ using \lemref{DGRS} 
to get
\begin{equation}\label{i21}
I_{21} \le 2^{10} \cdot 10^{13} \cdot \log 8 \cdot h_{\K} \pi^*(Q) \log Q 
+
\frac12 \pi^*(Q) \log Q .
\end{equation}
Let us now consider $I_{22}$.
We have
\begin{eqnarray*}
I_{22} 
& \leq &
6 \cdot 10^{30} h_{\K}Q 
\sum_{\delta \in \O_{\K} \atop   
1 \le \N(\delta) \leq 8Q^{3/2}}
\frac{1}{\N(\delta)Q^{1/4}}
\int_{\frac{Q\N(\delta)}{8}}^{Q\N(\delta)}
\frac{ \log u } {u}du \\
& \le &
3 \cdot \log 8 \cdot 10^{30} h_{\K}Q 
\sum_{\delta \in \O_{\K} \atop   
1 \le \N(\delta) \leq 8Q^{3/2}}
\frac{1}{\N(\delta)Q^{1/4}}
\log \frac{Q^2 \N(\delta)^2}{8} \\
& \le &
15 \cdot \log 8 \cdot 10^{30} h_{\K}\frac{Q \log (2Q)}{Q^{1/4}} 
\sum_{\delta \in \O_{\K} \atop   
1 \le \N(\delta) \leq 8Q^{3/2}}
\frac{1}{\N(\delta)}.\\
\end{eqnarray*}
Since $\log (2Q) \le 20 Q^{1/8}$, we have
$$
I_{22} \le 
3 \cdot \log 8 \cdot 10^{32} h_{\K}\frac{Q }{Q^{1/8}}
\sum_{\delta \in \O_{\K} \atop   
1 \le \N(\delta) \leq 8Q^{3/2}}
\frac{1}{\N(\delta)}.
$$
Now, using \eqref{sumndelta}, we get
\begin{eqnarray*}
I_{22} & \le & 10^{33} h_{\K}\frac{Q}{Q^{1/8}} 
(2\rho_{\K} \log 8Q + 2 \cdot \sqrt{3} \cdot 10^{15} h_{\K})\\
& \le &
10^{35} \rho_{\K} h_{\K}\frac{Q}{Q^{1/16}} 
 + \cdot 10^{49} h_{\K}^2  \frac{Q}{Q^{1/8}}.\\
\end{eqnarray*}
By \lemref{DGRS}, \lemref{comparisonterm},
the bound on $h_{\K}$ and the fact that
$Q \ge 8\exp(8\cdot 10^{45} |d_{\K}|)$  , we get
\begin{equation}\label{i22}
I_{22} 
\le
\frac12 \pi^*(Q) \log Q.
\end{equation}
Combining \eqref{i21} and \eqref{i22},
we get
\begin{equation}\label{i2}
I_2 \le 
(2^{12} \cdot 10^{13} \cdot h_{\K} + 1) \pi^*(Q) \log Q.
\end{equation}
Finally from \eqref{i1} and \eqref{i2},
we get 
\begin{eqnarray}\label{smalluprime}
	2\int_1^{Q^{5/2}} \left( \sum_{\P \subset \O_\K \atop{1 \leq  \N\mathfrak{\P} \leq u \atop{ \P \text{ principal}}}}\Lambda(\mathfrak{a})\sum_{\substack{\q \in \mathcal{T}_{\a} \atop \qq}}|H_{\q}(\K)| \right) \frac{du}{u^2}
\leq (2^{13} \cdot 10^{13} \cdot h_{\K} + 2) \pi^*(Q) \log Q.
\end{eqnarray}
\textbf{Contribution from $\psi(u)$:\;}	Using \thmref{GM2}, under the Riemann hypothesis we have  
for $u\geq 3$, 
$$\psi(u) \leq u +7\log{|d_\K|}\sqrt{u}\log u+\sqrt{u}\log^2{u}+19\sqrt{u}.
$$
Putting this in \eqref{1.2}, we get 
\begin{eqnarray}\label{smallpsi}
 2\int_1^{Q^{5/2}}\frac{\sum_{\qq}\psi(u)}{u^2}du\leq \left(5 \log Q +56\log|d_\K|+ 109 \right)\pi^*(Q)
\le
6 \cdot \pi^*(Q)\log Q.
\end{eqnarray}

We will now consider the case $u >Q^{5/2}$. Since 
	\begin{eqnarray*}
		\psi(u,\q,[1]_{\q})&=\displaystyle \sum_{\a \subset \O_{\K} \atop{1 \leq \N\a \leq u \atop \a \in [1]_{\q}}}\Lambda(\a) 
		=\displaystyle \sum_{\a \subset \O_{\K} \atop{1 \leq \N\a \leq u \atop{ (\a,\q)=\O_\K \atop (\a,\K(\q)/\K)=\sigma_0}}}\Lambda(\a) 
		=\displaystyle \sum_{\a \subset \O_{\K} \atop{1 \leq \N\a \leq u \atop (\a,\q)=\O_\K }}\mathbf{1}_{\sigma_o}((\a,\K(\q)/\K))\Lambda(\a),
	\end{eqnarray*} 
	where $(\a,\K(\q)/\K)$ denotes the Artin symbol of $\a$ and $\sigma_o$ the trivial element of $\operatorname{Gal}(\K(\q)/\K)$ with $\mathbf{1}_{\sigma_o}$ as its characteristic function.
Using \thmref{explicitchebotarev},
\lemref{size}, \lemref{Disc} and equation \eqref{CNF}, we get 
	\begin{align*}
\Big\lvert|H_{\q}(\K)|\psi(u,\q,[1]_{\q})-u\Big| &\leq \sqrt{u}\left(\left(\frac{\log u}{2\pi}+2 \right)\log |d_{\K(\q)}|+ \left(\frac{\log^2u}{4\pi}+4\right)|H_{\q}(\K)|\right)\\
&\leq 2\rho_{\K}\sqrt{|d_{\K}|}\N\q\Big(\log(\N\q|d_{\K}|)\log{u}+\log^2u\Big)\sqrt{u}.
\end{align*}
Therefore, for $Q\geq 8\exp(8\cdot 10^{45}|d_\K| )$,
\begin{eqnarray*}
&& 2  \int_{Q^{5/2}}^x  \frac{\sum_{\qq}\Big\lvert|H_{\q}(\K)|\psi(u,\q,[1]_{\q})- u\Big\rvert}{u^2}du \\
&\leq &  4\rho_{\K}\sqrt{|d_{\K}|} Q\pi^*(Q)\int_{Q^{5/2}}^x\frac{\log(Q|d_{\K}|)\log{u}+\log^2u}{u^{3/2}}du\\
&\leq & \frac12 \pi^*(Q)\log Q,
\end{eqnarray*}
where in the last step we have used 
\lemref{DGRS}. By \thmref{GM2}, for $Q\geq3$, we get 
$$\left|\psi(u)-u \right| \leq 7(\log{|d_\K|})\sqrt{u}\log u+\sqrt{u}\log^2{u}+19\sqrt{u}.$$
This implies that for $Q\geq 8\exp(8\cdot 10^{45}|d_\K| )$,
$$2\int_{Q^{5/2}}^x\frac{\sum_{\qq}|\psi(u)- u|}{u^2}du\leq (56\log|d_\K|+108)\pi^*(Q) \le\frac12 \pi^*(Q)\log{Q},
$$
where we have used equation \eqref{INT}.
Thus, the contribution of $Q^{5/2} <u<x$ to \eqref{1.2} is 
\begin{align}\label{largeu}
2\int_{Q^{5/2}}^x\frac{\sum_{\qq}\Big\lvert|H_{\q}(\K)|\psi(u,\q,[1]_{\q})- \psi(u)\Big\rvert}{u^2}du
& \leq \pi^*(Q)\log Q.
\end{align}
Combining \lemref{qbound}, \eqref{smalluprimepower}, \eqref{smalluprime}, \eqref{smallpsi} and \eqref{largeu}, we get 
\begin{align*}
\sum_{\qq} \Big\lvert \sum_{\chi \neq \chi_0}\Phi_\chi(x)\Big\rvert
 \leq (2^{13} \cdot 10^{13} h_{\K} + 6000h_{\K}^2 + 10) \pi^*(Q)\log Q.
\end{align*}
\end{proof}

Let $\q$ be a non-zero principal prime ideal of $\O_\K$.
We note that given a generalized Dirichlet character $\chi$ modulo $\q$ it is either primitive
or induced by a character of the class group (also called a generalized Dirichlet character modulo $\O_{\K}$).
We now consider the sum
\begin{eqnarray*}
	\sum_{\chi \neq \chi_0 } (\Phi_{\chi^*}(x) - \Phi_{\chi} (x)) 
	& = &
	\sum_{ \chi \bmod \q \atop{\chi^* \bmod \q }} (\Phi_{\chi^*}(x) - \Phi_{\chi} (x)) 
	+
	\sum_{ \chi \bmod \q \atop{ \chi^* \bmod \O_{\K} \atop{\chi \neq \chi_0 }}} (\Phi_{\chi^*}(x) - \Phi_{\chi} (x)) \\
	& = &
	\sum_{ \chi \bmod \q \atop{\chi^* \bmod \O_{\K} \atop{\chi \neq \chi_0 }}} (\Phi_{\chi^*}(x) - \Phi_{\chi} (x))\\
	& = &
	\sum_{ \chi^* \bmod \O_{\K} \atop{\chi \neq \chi_0}} \frac{1}{x-1} \int_1^x \left( \sum_{\N\a \le t, \atop{ \q \mid \a}} \frac{\Lambda(\a)}{\N\a} \chi^*(\a) \right) dt.
\end{eqnarray*}
Since $\q$ is a prime ideal, it follows that the inner sum only consists of terms for which $\a = \q^k$. Further since $\q$ is principal so is $\a$.
Therefore, we get
\begin{eqnarray*}
	\sum_{\chi \neq \chi_0} (\Phi_{\chi^*}(x) - \Phi_{\chi} (x))
	& = &
	\frac{1}{x-1} \int_1^x \left( \sum_{\N\a \le t, \atop{ \q \mid \a}} \frac{\Lambda(\a)}{\N\a} \sum_{ \chi^* \bmod \O_{\K}  \atop{\chi \neq \chi_0}} \chi^*(\a) \right) dt \\
	& \le &
	\frac{h_{\K}}{x-1} \int_1^x \left( \sum_{k \ge 1 \atop{ \N\q^k \le t}} \frac{ \log \N\q}{\N\q^k} \right) dt\\
	& \le & \frac{h_{\K} x}{x-1}  \cdot \frac{ \log \N\q}{\N\q -1 }. 
\end{eqnarray*}
Therefore for $x \ge 2$, using  \lemref{qbound}, we obtain
\begin{equation}\label{phidiff}
\sum_{\frac{Q}{2} < \N\q \le Q} \Bigg|\sum_{\chi \neq \chi_0 } (\Phi_{\chi}(x) - \Phi_{\chi^*} (x)) \Bigg| \le  2h_{\K} \left( 14 + 4\frac{e^{75}|d_\K|^{1/3}(\log|d_\K|)^2}{\rho_\K} \right).
\end{equation}
\smallskip
\begin{prop}\label{Phi2}
	Let $Q \geq 2 $ and $x \geq Q^4$. Let $\chi$ be a generalized Dirichlet character modulo $\q$, where $\q$ is a non-zero principal prime ideal in $\O_\K$. Also let
	$\frac{L'}{L}(1,\chi)$ denote the logarithmic derivative of the corresponding Hecke L-function evaluated at $1$.
	 Then under GRH, we have
	$$\sideset{}{'}\sum_{\qq}\sum_{\chi \neq \chi_0}\Bigg\lvert \frac{L'}{L}(1,\chi^*)+\Phi_{\chi^*}(x)\Bigg\rvert 
	< 2010 \rho_\K \sqrt{|d_\K|} \frac{ \pi^*(Q) \log( 5|d_\K|Q)}{Q}
	,$$
	where the outer sum is over principal prime ideals in the interval $(\frac{1}{2}Q, Q] $ and $ \pi^*(Q)$ denotes the number of principal prime ideals with norm in the interval $(\frac{1}{2}Q, Q]$.
\end{prop}

\begin{proof}
	Using \thmref{ims}, we have	
	\begin{equation}\label{1.5}
		\frac{L'}{L}(1, \chi^*) + \Phi_{\chi^*}(x)= \frac{1}{x-1}\sum_{\rho_{\chi^*}}\frac{x^{\rho_{\chi^*}}-1}{\rho_{\chi^*}(1-\rho_{\chi^*})}+\log{\frac{x}{x-1}}+\frac{1}{x-1}\log{x},
	\end{equation}
	where the sum is over all non trivial zeros of $L(s,\chi^*)$ counted with multiplicities. 
	Using \propref{zeroes}, we get
	$$N_{t,\chi^*}=\#\{\rho=\frac{1}{2}+i \beta : | \beta -t|\leq 1, L(\rho, \chi^*)=0\} < 100 \log (|d_\K| \N\q~(|t|+ 2)).$$ 
For a non-trivial zero $\rho_{\chi^*}=\frac{1}{2}+i \beta$ of $L(s, \chi^*)$, we have $|\rho_{\chi^*}(1- \rho_{\chi^*})|= \left(\frac{1}{2}\right)^2+  \beta^2$. For $|\beta| \leq 1$, we get 
\begin{equation}\label{green}
		\sum_{\rho_{\chi^*}=1/2+i\beta \atop |\beta|\leq 1}\frac{1}{|\rho_{\chi^*}(1- \rho_{\chi^*})|}\leq 4 N_{0,\chi^*}  < 400 \log (2|d_\K| \N\q).
	\end{equation}
	Also, for $|\beta|>1$,
	 we get
	\begin{align*}
		\sum_{\rho_{\chi^*}=1/2+i\beta \atop |\beta|> 1}\frac{1}{|\rho_{\chi^*}(1- \rho_{\chi^*})|} & \leq \displaystyle \sum_{k=1}^{\infty}\sum_{\rho_{\chi^*}=1/2+i\beta \atop k<|\beta|\leq k+1}\frac{1}{\beta^2}  
		\leq \sum_{k=1}^{\infty} 
		\frac{1}{k^2} ( N_{k,\chi^*} + N_{-k,\chi^*}) \nonumber \\
		& < 200 \displaystyle \sum_{k=1}^{\infty} 
		\frac{\log (|d_\K| \N\q~(k+ 2))}{k^2}  .
	\end{align*}
We have
\begin{eqnarray*}
	\sum_{k=1}^{\infty} 
	\frac{\log (|d_\K| \N\q~(k+ 2))}{k^2} 
	& \le &
	\zeta(2)\log (|d_\K| \N\q) 
	+
	\sum_{k=1}^{\infty} \frac{\log(k+2)}{k^2} 
	 \le 
	3 \log(5|d_{\K}| \N\q).
\end{eqnarray*}

It implies that
	\begin{align}\label{gre}
	\sum_{\rho_{\chi^*}=1/2+i\beta \atop |\beta|> 1}\frac{1}{|\rho_{\chi^*}(1- \rho_{\chi^*})|}
	\leq 600 \log (5|d_\K| \N\q).
	\end{align}
Plugging \eqref{green} and \eqref{gre} in \eqref{1.5}, we get
	\begin{align*}
		\Big|\frac{L'}{L}(1, \chi^*) + \Phi_{\chi^*}(x)\Big| 
	&< \frac{\sqrt{x}+1}{x-1}\sum_{\rho_{\chi^*}}\frac{1}{|\rho_{\chi^*}(1-\rho_{\chi^*})|} + \log{\frac{x}{x-1}}+  \frac{2\log x}{x} \\
	& <  \frac{2000}{\sqrt{x}} \log(5|d_\K| \N\q)+ \frac{1}{x-1} +  \frac{4}{\sqrt{x}}\\
	& < \frac{2010}{\sqrt{x}} \log(5|d_\K| \N\q).
\end{align*}
Summing over $\chi \neq \chi_0$ modulo $\q$ and then over the principal prime ideals $\q $ with $\N\q \in (\frac{1}{2}Q,Q]$, we get 
		\begin{align*}
		\sideset{}{'}\sum_{\qq}\sum_{\chi \neq \chi_0}	\Big|\frac{L'}{L}(1, \chi^*) + \Phi_{\chi^*}(x)\Big| 
		&<  \frac{2010}{\sqrt{x}} h_\K Q\pi^*(Q) \log( 5|d_\K|Q)\\
		&\leq \frac{2010}{\sqrt{x}}\rho_\K\sqrt{|d_\K|}Q\pi^*(Q) \log( 5|d_\K|Q),
	\end{align*}
	where we use equation \eqref{CNF} in the last step. Choosing $x \geq Q^4$, we get our desired result.
\end{proof}

\smallskip

\smallskip
\subsection{Proof of \thmref{JR}}
Let $\chi$ be a generalized Dirichlet character modulo $\q$, where $\q$ is a non-zero principal prime ideal in $\O_\K$. 
We have (see Chapter VIII, \cite{CF} for more details)
 \begin{equation} \label{lprod}
\zeta_{\K(\q)}(s)= \prod_{\chi \mod \q}L(s, \chi^*), 
\end{equation}
where the product is over all generalized Dirichlet characters modulo $\q$. 
For any number field $\K$, an equivalent definition of $\gamma_{\K}$
 is given by 
$$
\gamma_\K = \lim_{s \rightarrow 1}\left(
\frac{\zeta'_{\K}}{\zeta_{\K}}(s) + \frac{1}{s-1} \right).
$$ 
Therefore, taking the logarithmic derivative
on both sides of \eqref{lprod}, we get
$$\gamma_{\K(\q)}=\gamma_{\K} + \sum_{\chi \neq \chi_0}\frac{L'}{L}(1,\chi^*).$$
It now follows that
	$$\gamma_{\K(\q)}=\gamma_{\K}-\sum_{\chi \neq \chi_0}\Phi_{\chi}(x)+\sum_{\chi \neq \chi_0}\Big(\Phi_{\chi}(x)-\Phi_{\chi^*}(x)\Big)+ \sum_{\chi \neq \chi_0}\Big(\frac{L'}{L}(1,\chi^*)+\Phi_{\chi^*}(x)\Big).$$
Therefore,
	$$|\gamma_{\K(\q)}|\leq |\gamma_{\K}| + \left|\sum_{\chi \neq \chi_0}\Phi_{\chi}(x)\right|+ \left|\sum_{\chi \neq \chi_0}\Big(\Phi_{\chi}(x)-\Phi_{\chi^*}(x)\Big)\right|+\left| \sum_{\chi \neq \chi_0}\Big(\frac{L'}{L}(1,\chi^*)+\Phi_{\chi^*}(x)\Big)\right|.$$
Using \thmref{Phi1}, \eqref{phidiff} \propref{Phi2} and the fact that $Q\geq 8\exp(8\cdot 10^{45}|d_\K| )$, we get
	\begin{eqnarray*}
		\sideset{}{'}\sum_{\qq}|\gamma_{\K(\q)}| & < & |\gamma_{\K}|\pi^*(Q)  + ( 6000h_{\K}^2+ 2^{13} \cdot 10^{13} h_{\K}  + 10) \pi^*(Q)\log Q \\
		& + &  2h_{\K} \log Q
		 +   \pi^*(Q) \log Q.
	\end{eqnarray*}
	This implies that
	\begin{eqnarray*}
 \frac{1}{\pi^*(Q)}	\sideset{}{'}\sum_{\qq}|\gamma_{\K(\q)}|
 & < &
 |\gamma_{\K}| +
(6000h_{\K}^2+ 10^{17}h_\K  + 11) \log Q.
	\end{eqnarray*}

\section{\large{Concluding Remarks}}
We conclude by showing that
infinitely many narrow ray class fields corresponding to principal prime moduli of an
imaginary quadratic field are non-abelian. Given any narrow ray class field $\K(\q)$
of $\K$ modulo an integral ideal $\q$, the Hilbert class field of $\K$ is a subfield of $\K(\q)$ 
containing $\K$. By a result of Ankeny and Chowla \cite{AC}, we know that there are infinitely
many imaginary quadratic fields whose class number is divisible by $3$. By genus theory (see chapter 1 of \cite{MI})
if the Hilbert class field of $\K$ is abelian over $\Q$, the class number must be a power of $2$. Therefore for the
family of fields suggested by the result of Ankeny and Chowla, the Hilbert class field must be
non-abelian and therefore all the narrow ray class fields of such an imaginary quadratic are also non-abelian.

\bigskip
\noindent{\bf Data Availability Statement:} Data sharing is not applicable to this article as no datasets were generated or analysed during the current study.

\bigskip
\noindent{\bf Conflict of Interest Statement:} The authors certify that there is no conflict of interest.

\bigskip
\noindent
{\bf Acknowledgments.} 
All the authors
would like to thank Prof. Sanoli Gun for suggesting the problem and comments on the manuscript. 
They would also like to thank Prof. Purusottam Rath for helpful comments on earlier versions of the draft. The first author would also like to thank Prof. Alessandro Languasco for drawing her attention to the reference \cite{dixit-errata}.
Part of this work was done during visits of the first author to the Sorbonne university
and the second author to the Sorbonne and the Paris Cit{\'e} universities under the
support of the Indo-French Program in Mathematics (IFPM). The first and second authors would like to
thank the IFPM for partial financial support. 
They would also like to thank the DAE number theory plan project.
All the authors thank the Institute of Mathematical Sciences for providing academic facilities. 
The third author would like to thank the Chennai Mathematical Institute and the Indian Statistical Institute for providing academic support.
Further, the third author would like to thank
 the National Board for Higher Mathematics and the Department of Science and Technology for financial support.
 Finally the authors would like to thank the referee for a careful reading of the paper and a number of relevant suggestions.

\end{document}